\documentclass[11pt]{article}

\usepackage[centertags]{amsmath}
\usepackage{amsthm}
\usepackage{amssymb}
\usepackage{amscd}
\usepackage{eucal}
\usepackage{epsfig}
\usepackage{verbatim}

\usepackage{latexsym,enumerate,amsxtra}
\usepackage{color}

\usepackage[margin = 1.1in]{geometry}

\newtheorem{thm}{Theorem}[section]
\newtheorem{theorem}[thm]{Theorem}

\newtheorem{lemma}[thm]{Lemma}
\newtheorem{proposition}[thm]{Proposition}
\newtheorem{definition}[thm]{Definition}

\newcommand{\bsigma}{\boldsymbol{\sigma}}

\newcommand{\balpha}{\boldsymbol{\alpha}}
\newcommand{\bQ}{\boldsymbol{Q}}
\newcommand{\bV}{\boldsymbol{V}}
\newcommand{\bW}{\boldsymbol{W}}

\newcommand{\bA}{\boldsymbol{A}}

\newcommand{\bu}{\boldsymbol{u}}
\newcommand{\bw}{\boldsymbol{w}}

\newcommand{\bm}{\boldsymbol{m}}
\newcommand{\blambda}{\boldsymbol{\lambda}}

\newcommand{\brho}{\boldsymbol{\rho}}
\newcommand{\bphi}{\boldsymbol{m}}

\newcommand{\K}{\mathcal K}

\newcommand{\I}{\mathcal I}

\newcommand{\cC}{\mathcal C}
\newcommand{\J}{\mathcal J}
\newcommand{\cP}{\mathcal P}
\newcommand{\cH}{\mathcal H}
\newcommand{\cF}{\mathcal F}

\newcommand{\Lloc}{\mathcal L_{\mathrm{local}}}
\newcommand{\Lglo}{\mathcal L_{\mathrm{global}}}

\newcommand{\Z}{\Lloc}

\begin{document}

\title{From Local to Global Stability in Stochastic Processing Networks through Quadratic Lyapunov Functions}

\author{A.\ B. Dieker\thanks{Antonius B. (Ton) Dieker is with the H.\ Milton Stewart School of Industrial and Systems Engineering, Georgia Institute of Technology,
Atlanta GA USA 30332. Email: ton.dieker@isye.gatech.edu.}
 \and J.\ Shin\thanks{Jinwoo Shin is with the Mathematical Sciences Department, IBM T.\ J.\ Watson Research,
Yorktown Heights, NY USA 10598. Email: jshin@us.ibm.com.}
}

\maketitle

\begin{abstract}

We construct a generic, simple, and efficient scheduling policy for stochastic
processing networks, and provide a general framework to establish its stability.
Our policy is randomized and prioritized: with high probability it prioritizes jobs
which have been least routed through the network. We show that the network is globally
stable under this policy if there exists an appropriate quadratic `local' Lyapunov function
that provides a negative drift with respect to nominal loads at servers.
Applying this generic framework, we obtain stability results for our policy in many important examples of stochastic processing
networks: open multiclass queueing networks, parallel server networks, networks of input-queued switches,
and a variety of wireless network models with interference constraints.
Our main novelty is the construction of an appropriate `global' Lyapunov function from
quadratic `local' Lyapunov functions, which we believe to be of broader interest.

\end{abstract}

\section{Introduction}
The past few decades have witnessed a surge in interest on the design and
analysis of scheduling policies for stochastic
networks, e.g., \cite{Kelly,Tassiulas,McKeown,Mo,harrison1,Dai08}.
One of the key insights from this body of work is that natural scheduling
policies can lead to instability even when
each server is nominally underloaded \cite{Kumar90,Lu91,RybkoStolyar}.
(There are several notions of stability for stochastic networks, but
they intuitively entail that, in some sense,
the number of customers in the system does not grow without bounds.)
This insight stimulated a search for
tools that can characterize the stability regions of scheduling policies, i.e.,
the exact conditions on the arrival and service rates under which a network is stabilized by a policy.

It is the objective of this paper to study a question of a different kind:
is it possible to construct a generic, simple, and efficient scheduling policy
for stochastic processing networks, which leads to a (globally) stable
network if all servers are (locally) nominally underloaded in some sense?
To our knowledge, we are the first to answer this question within the setting of stochastic processing networks,
which constitute a large class of stochastic networks capable of
modeling a variety of networked systems for communication, manufacturing, and service systems
(e.g., \cite{harrison1}).
To investigate this question, the key is to determine if and how jobs
from different parts of the network should be treated differently when they share
the same buffer.

Various existing scheduling policies are `throughput optimal'
in the sense that they achieve the largest possible stability region,
but these policies suffer from significant drawbacks.
They typically obtain the desired stability by framing the contention resolution between
buffers as an appropriate global optimization problem.
This optimization problem requires central coordination between network entities,
and it is computationally hard to solve if the network is large.
The resulting policies, such as the max-weight policy \cite{McKeown,shah2010switched} and the back-pressure policy
\cite{Tassiulas,Dai08}, are not scalable and cannot cope efficiently with large networks.
As a result, these throughput optimal policies do not provide a satisfactory answer to the aforementioned question.

The computational challenges surrounding existing throughput optimal policies
motivate the search for easily implementable scheduling policies with provable
performance guarantees but not necessarily with the throughput optimality property.
This has led to the analysis of simple greedy scheduling policies for a variety of special
classes of stochastic processing networks: open multiclass queueing networks \cite{Dai95},
input-queued switches \cite{Dai00}, and wireless network models \cite{Srikant09}.
These policies play an important role in this paper, since our most critical assumption roughly requires any
local network component to be nominally underloaded under
any `maximal' greedy scheduling policy.

Our main contribution is a randomized scheduling policy
for stochastic processing networks which only requires
coordination within local components (e.g., service stations), and which
is computationally attractive since it is a kind of priority policy.
With high probability, our policy prioritizes jobs which have been least routed, and
we therefore call our policy the $\varepsilon$-Least Routed First Served
($\varepsilon$-LRFS) policy. Here $\varepsilon\ge 0$ is a small number
which helps to make the meaning of ``high probability'' precise.

Our main technical tool is a novel framework to
construct a `global' Lyapunov function for a stochastic processing network through
appropriate `local' Lyapunov functions.
If the local Lyapunov functions yield stability
of the corresponding `local' network components,
then the global Lyapunov function allows
us to conclude that the whole network is stable.
A critical feature of our framework is that the Lyapunov functions we work with are quadratic.
Through examples,
we show that quadratic local Lyapunov functions can readily be found for wide classes of networks.
We refer to \cite{diekergao:OU,Kumar95} for other uses of quadratic Lyapunov functions.

Our approach to construct an appropriate `global' Lyapunov function
using `local' quadratic Lyapunov functions
contrasts with the popular fluid model methodology for establishing stability of
stochastic networks \cite{Dai95,dupuiswilliams,RybkoStolyar}.
The fluid model framework essentially reduces the question of stochastic stability to a question
of a related deterministic (fluid) system.
In the case of reentrant lines, our policy reduces to the First Buffer
First Served (FBFS) policy, which has been proven to be stable via `inductive' fluid arguments \cite{Dai96}.
A similar fluid induction argument can be expected to work for general multiclass networks (modulo
some technical arguments),
but a fluid induction argument cannot be expected to work in general.
A disadvantage compared to fluid models is that we have to keep track of
detailed system behavior such as remaining service times,
but therein also lies the power of our approach.
The generality of our framework presents challenges to the use of fluid methods,
and instead we work directly with a global Lyapunov function.
The connection with fluid techniques
is discussed in more detail in Section~\ref{sec:connectionfluid}.

Although we believe that our methodology could be of much wider use,
we have chosen to work out two special classes of stochastic processing networks
in order to describe the implications of our techniques in relatively simple yet powerful settings:
parallel server networks (including multiclass queueing networks)
and communication network models (including networks of input-queued switches and wireless networks).

Our work is related to a paper by Bramson \cite{Bramson00}, who shows
that open multiclass queueing networks can be stabilized
by the Earliest Due Date First Served (EDDFS) policy
under the `local' condition that every processing unit is nominally underloaded.
However, it is not known whether the EDDFS policy achieves similar results
beyond the open multiclass queueing network.
In fact, it is {\sl a priori} unclear how to formulate a `local' condition
for general stochastic processing networks;
our notion of `local' Lyapunov function plays this role in the present paper.

This paper is organized as follows.
Section~\ref{sec:defSPN} introduces the class of stochastic processing networks we
study in this paper. Section~\ref{sec:mainresult} formally introduces the $\varepsilon$-LRFS
policy and presents our main results.
Section \ref{sec:overview} describes our main idea
using a simple network, the Rybko-Stolyar network \cite{RybkoStolyar}.
We specialize our result to parallel server and communication network models in Section~\ref{sec:application}.
All proofs are given in Section~\ref{sec:pfthmmain}.

\section{A Class of Stochastic Processing Networks}\label{sec:defSPN}

A stochastic processing network (SPN) consists of a set
$\I=\{1,\dots,I\}$ of buffers, a set $\J=\{1,\dots,J\}$ of
activities and a set $\K=\{1,\dots,K\}$ of processors. Each buffer
has infinite capacity and holds jobs that await service. The SPNs we
study in this paper have the feature that activity $j\in \J$ can
only process jobs from a single buffer $i_j\in\I$, and that $j$
requires simultaneous possession of a set $\K_j\subseteq \K$ of
processors. Let $\J_i$ be the set of activities capable of processing
buffer $i$, i.e.,
$$\J_i~=~\{j\in\J:i_j=i\}.$$
We say that two buffers $i$ and $\ell$ are {\em
activity-interchangeable}
if $\{\K_j: j\in \J_i\} =\{\K_j: j\in
\J_\ell\}$. We also say that buffers $i$ and $\ell$
are {\em processor-independent} if $\bigcup_{j\in \J_i} \K_j$ and
$\bigcup_{j\in \J_\ell} \K_j$ are disjoint.

\paragraph{Network state.}
We let $Q_i(t)\in \mathbb Z_+$ be the queue length of buffer $i$ at time $t$, i.e.,
the number of jobs waiting in buffer $i$ excluding those being
processed. We write $V_i(t)\in \mathbb R_+$ for the sum of the remaining service
requirements over all jobs in buffer $i$ which are currently being
processed at time $t$. We use $\sigma_j(t)$ to denote the activity level
of activity $j$ at time $t$, where
we assume that $\bsigma(t)=[\sigma_j(t)]\in\{0,1\}^J$, meaning that
each activity is either fully employed or not employed at all. Thus,
we say that
the network is {\em non-processor-splitting}.

\paragraph{Routing.}
After departing from a buffer $i$, a job joins buffer $\ell\in\I$
with probability $P_{i\ell}$ and departs from the network with
probability $1-\sum_{\ell} P_{i\ell}$ (independently of everything
else). We write $P$ for the $I\times I$ matrix of routing probabilities.

\paragraph{Resource allocation.}
Each activity $j$ decreases the remaining service requirement of the
job it is processing at rate $\beta_j>0$ if $\sigma_j(t)=1$. It is
not allowed for an activity to be interrupted before it finishes the service
requirement of the job it is working on, i.e., the network is {\em
non-preemptive}.
We do not allow for multiple activities to work on the same job simultaneously.
Furthermore, we assume that each
processor $k$ has unit capacity, i.e.,
\begin{equation*}
\sum_{j\in\J : k\in \K_j} \sigma_j(t)~\leq~ 1.
\end{equation*}
We note that this unit capacity assumption is not restrictive since
in our non-processor-splitting network of $\bsigma(t)\in\{0,1\}^I$,
a single processor $k$ with capacity $c_k\in \mathbb{N}$ can be
replaced by $c_k$ copies with unit capacity, and identical activity
structure inherited from the original processor. We further define
the service rate vector $s(\bu)=\left[s_i(\bu)\right]\in
\mathbb{R}_+^I$ for a (scheduling) vector $\bu=[u_j]\in\{0,1\}^J$
through
$$s_i(\bu)~:=~\sum_{j\in \J_i}u_j\beta_j.$$

\paragraph{External arrivals.}
There are external arrivals to at least one buffer. For $0\le s\le
t$, let $A_i(s,t)$ be the number of external jobs arriving at buffer $i$
during the time interval $[s,t)$.
We assume that $E\left[A_i(0,t)^2\right]<\infty$ for all $t<\infty$
(i.e., bounded second moment)
and
$$\limsup_{t\to\infty}E\left[A_i(t,t+1)~\Big|~\big\{A_1(0,s),\dots,A_I(0,s):0\leq s\leq t\big\}\right]<\alpha_i ,
$$for constant
$\alpha_i\in(0,\infty)$.
We note that we allow for
dependencies in the random processes $A_i(\cdot)$ and $A_\ell
(\cdot)$ for two buffers $i$ and $\ell$. Even though it is possible
for the $\alpha_i$ to exceed the external arrival rates, it is
convenient to interpret $\alpha_i$ as the external arrival rate at
buffer $i$. Similarly abusing terminology, we let
$\blambda=[\lambda_i]$ be the {\em effective} arrival rate vector,
i.e.,
$$\blambda~:=~\left(I+P+P^2+\cdots\right)\balpha,$$
where $\balpha=[\alpha_i]$. We also say that all routes are bounded
(in length) if $\inf\left\{d\in\mathbb{N}:P^d=0\right\}<\infty.$

\paragraph{Service requirements.}
Once a job in buffer $i$ is selected for processing by activity
$j\in \J_i$, it requires service for a random amount of time.
We assume that all service times are independent, and that they
are independent of the routing and external arrival processes.
We also suppose that
the service time distribution only depends on the buffer from which the job is
processed.
Writing $\Gamma_i$ for a generic processing time at buffer
$i$, we assume that
$$
E[\Gamma_i]=m_i\qquad\mbox{and}\qquad E\left[\Gamma_i^2\right] <\infty,$$for
constants $m_i\in(0,\infty)$. Let $\brho=[\rho_i]$ denote the
nominal load, i.e., $\rho_i=\lambda_i m_i$. We write $W_i(t)$
for the (expected) immediate workload in buffer $i$, which we define to be
\begin{equation}\label{eq:defW}
W_i(t)~=~m_iQ_i(t)+V_i(t). \end{equation}
This is the expected amount of work in the $i$-th buffer given $Q_i(t)$ and $V_i(t)$.
We note that our notion of immediate workload is the conditional expectation of a more common
definition, and that immediate workload is defined for each buffer (as opposed to resource).
We say that the network is {\em synchronized} if for all $s,t\geq 0$, $i\in
\I$ and $j\in\J$,
$$A_i(s,t)=A_i(s,\lfloor t\rfloor),\qquad\mbox{ and }\qquad
\Gamma_i=m_i=\beta_j=1\quad\mbox{with probability $1$}.$$
That is, in a synchronized network,
arrivals and service completions only occur at integer time epochs. 

\paragraph{Maximal scheduling policies.}

We say that activity $j$ is {\em maximal} in a (scheduling) vector
$\bu=[u_j]\in \{0,1\}^J$
if there exists a processor $k\in \K_j$ such that
\[
\sum_{\ell\in\J:\, k\in \K_\ell} u_\ell~=~ 1,
\]
i.e., either activity $j$ uses each of the processors in $\K_j$ under
schedule $\bu$ or it cannot be employed without violating the unit
capacity constraint for some processor in $\K_j$.
Given a non-negative vector $\bw=[w_i]\in \mathbb{R}_+^I$, activity
$j$ is called {\em maximal} in $\bu$ with respect to $\bw$ if
$w_{i_j}=0$ or $j$ is maximal in $\bu$.
An activity is non-maximal if it is not maximal.
If all activities are
maximal in $\bu$ with respect to $\bw$, we simply say that $\bu$
is maximal with respect 
to $\bw$.
We write $\mathcal M(\bw)$ for the set of maximal scheduling vectors
with respect to $\bw$, i.e., $$\mathcal
M(\bw)~:=~\left\{\bu\in\{0,1\}^J:\mbox{$\bu$ is maximal with respect
to $\bw$}\right\}.$$
 Finally, we say that a scheduling policy is
maximal if, under the policy, $\bsigma(t)\in \mathcal M(\bQ(t))$ for
all $t\geq 0$.

\section{Main Result}\label{sec:mainresult}
This section describes the scheduling policies which play a central
role in this paper, and presents our main stability result.

\paragraph{Scheduling policies.}
Our policies require that each job maintains a `counter' for the
number of times it has been routed so far, where we follow the
convention that counters start from $1$. The counter of a job is
increased even when a job is routed to a buffer it has previously visited,
 so the counter of a job could differ from the number of different
stations it has visited.
We also consider a partition $\{\I^{(h)}: h\in\mathcal H=\{1,\dots, H\}\}$
of buffers into components,
such that $$\I=\bigcup_{h\in \cH}\I^{(h)},$$
where any two buffers from different components are processor-independent.
One possible choice for the partition is $\I=\I^{(1)}$, but, as becomes apparent from the 
description of our policies below, a finer partition makes our policies more `distributed'.
It is important to note that we allow for routing between components.
We assume that each component $\I^{(h)}$ for $h\in \cH$
maintains a `timer' $\mathcal T^{(h)}(t)\in[0,1]$ at time $t$
which decreases at unit rate if $\mathcal T^{(h)}(t)>0$.

To describe our maximal scheduling policies, we need the following notation.
We write $\bQ^{(h)}(t)=\left[Q_i(t)I_{\I^{(h)}}(i)\right]$
for the queue length information in component $\I^{(h)}$, where
$I_S\in\{0,1\}$ is the indicator function of the set $S$, i.e.,
$$I_S(i)~=~\begin{cases}1&\mbox{if $i\in S$}\\
0&\mbox{otherwise}\end{cases}.$$
The following policy plays a key role throughout this paper.

\begin{definition}[{\bf Least Routed First Served (LRFS) policy}]
For each $h\in\mathcal H$, whenever a new arrival and/or service completion occurs at time $t$ in component $\I^{(h)}$,
execute the following algorithm immediately after all arrivals and service completions have occurred:
\begin{itemize}
\item[1.] Find the set $\Sigma$ of non-maximal activities in $\bsigma(t)$
with respect to
$\bQ^{(h)}(t)$.
\item[2.] Find a job with the smallest counter among those in buffers $\{i_j: j\in\Sigma\}$, where ties are broken arbitrarily.
\item[3.] Choose an arbitrary activity $j\in \Sigma$ to process the job identified in step 2, i.e., set $\sigma_{j}(t)=1$.
\item[4.] Repeat steps 1--3 until $\bsigma(t)$
is maximal with respect to $\bQ^{(h)}(t)$.
\end{itemize}
\end{definition}
\noindent We next introduce the maximal scheduling policy which is of primary interest in this paper.
It uses a small parameter $\epsilon>0$ in order to deal with
unbounded route lengths.

\begin{definition}[\bf $\varepsilon$-Least Routed First Served ($\varepsilon$-LRFS) policy]\label{def:epLRFS}
For each $h\in\mathcal H$, whenever a new arrival and/or service completion occurs at time $t$ in component $\I^{(h)}$,
execute the following algorithm immediately after all arrivals and service completions have occurred:
\begin{itemize}
	\item[1.] Find a buffer $i\in \I^{(h)}$ containing a job with the largest counter among those in buffers $\I^{(h)}$.
\item[2.] Find the set of non-maximal activities $\Sigma\subset \J_i$ in $\bsigma(t)$. 
\item[3.] If $\Sigma\neq \emptyset$ and $\mathcal T^{(h)}(t)=0$,
\begin{itemize}
	\item[3-1.] With probability $\varepsilon$, choose an arbitrary activity $j\in \Sigma$ to process
                    the job identified in step 1, i.e., set $\sigma_j(t)=1$. With probability $1-\varepsilon$, do nothing.
	\item[3-2.] Set $\mathcal T^{(h)}(t)=1$.
\end{itemize}
\item[4.] Execute the LRFS policy for component $\I^{(h)}$. 
\end{itemize}
\end{definition}
\noindent We remark that the $\varepsilon$-LRFS policy
is identical to the LRFS policy when either $\varepsilon=0$
or $\mathcal T^{(h)}(t)>0$,
i.e., execution of the first three steps is not necessary in these cases.
Since we assume that each timer decreases at unit rate,
step 3-1 can be executed at most once per component in any time
interval of unit length.

\paragraph{The network process.}
Write $Q_{i,c}(t)$ for the queue length of jobs with counter $c$ in
buffer $i$ at time $t$. Let $V_{i,c}^j(t)$ be the remaining service
requirement of the job with counter $c$ in buffer $i$ at time $t$ if
it is processed by activity $j$, and set $V_{i,c}^j(t)=0$ if $j$ is
not processing a job with counter $c$.
The network state is described by
$$X(t)~=~\left[Q_{i,c}(t), V_{i,c}^j(t), \mathcal T^{(h)}(t)~:~i\in \I, j\in
\J, c\in \mathbb N,h\in \cH\right]~\in~\Omega_X:=\left(\mathbb{Z}_+^I\times
\mathbb{R}_+^{I\times J}\right)^{\infty}\times [0,1]^H.$$
Note that $X(t)$ does not encode information
on the external arrival processes. In particular, $\{X(t)\}$ is non-Markovian in general.
We impose the convention that $\{X(t)\}$ has right-continuous sample paths.
We define a norm on
$\Omega_X$ through
\begin{eqnarray*} |X(t)|~=~\sum_{i\in\I, c\in\mathbb{N}}
Q_{i,c}(t)+\sum_{i\in\I, c\in\mathbb{N},j\in\J}
V_{i,c}^j(t).
\end{eqnarray*}
We assume that $X(0)\in \Omega_X^*$ where
$\Omega_X^*=\{X\in\Omega_X:|X|<\infty\}$. In synchronized networks,
one has $V^j_{i,c}(t)\in\{0,1\}$ for all $t\in\mathbb{Z}_+$ under
the $\varepsilon$-LRFS policy. We refer to the process
$\{X(t):t\in\mathbb{R}_+\}$ operating under the $\varepsilon$-LRFS
policy as the $\varepsilon$-LRFS process.
If $\varepsilon=0$, then we simply refer to this process as the LRFS process.

\paragraph{Network stability.}
This paper uses the following notion of stability.
\begin{definition}
The $\varepsilon$-LRFS process $\{X(t)\}$ is called
queue-length-stable if
\begin{equation}\label{eq:queuestable}
	\limsup_{t\to\infty}\frac1t\int_0^t
E\left[|X(s)|\right]ds<\infty,\end{equation} for any given initial state
$X(0)\in\Omega_X^*$.

\end{definition}

We establish the queue-length-stability by constructing
appropriate Lyapunov functions.
Under some additional assumptions
on the arrival processes and service time distributions, these Lyapunov functions can also be used
to establish positive recurrence of
the $\varepsilon$-LRFS process, cf.~Condition (A3) in \cite{dai1995stability}.
It is outside of the main scope of the current paper to work out the details.
We also remark that our proof can find an explicit
finite constant for the right-hand side of \eqref{eq:queuestable},
but this requires tedious bookkeeping and we
therefore do not carry out this analysis.

To show the desired stability, we need the following
notion of a `local' Lyapunov function.
\begin{definition}\label{def:local}
We say that $\Lloc:\mathbb{R}_+^I\to\mathbb{R}_+$ is a {\em local
Lyapunov function} with slack parameter $\varepsilon\geq 0$ if there
exist constants $\mathcal C\geq 0$, $\eta>0$ such that for every
pair $(\bw,\bu)\in\mathbb{R}_+^I\times \{0,1\}^J$ satisfying
$\bu\in\mathcal M(\bw)$,
\begin{equation}
\label{eq:localcondition} \mathcal L_{\text{\em
local}}\Big(\bw+\brho+\varepsilon \bm-s(\bu)\Big)~ \leq~\mathcal
L_{\text{\em local}}(\bw)-\eta \|\bw\|_1+\mathcal
C,
\end{equation}
where $\bm=[m_i]\in\mathbb R^{I}$.
\end{definition}
The reason for this nomenclature is that in a `local' network, i.e.,
a network without routing ($P=0$), the above inequality for $\Lloc$
provides the desired negative drift condition in the Foster-Lyapunov criteria \cite{foss2004},
which implies network stability.
Here $w_i$, $\rho_i$, and $s_i(\bu)$ can
be interpreted as the immediate workload, the 
external workload arrival rate and the workload processing rate at buffer $i$,
respectively. We refer to Sections~\ref{sec:overview} and
\ref{sec:application} for examples of local Lyapunov
functions. 

Now we are ready to state the main theorem of this paper, which
establishes `global' stability using a `local' quadratic Lyapunov
function.
\begin{theorem}
\label{thm:main} Suppose that there exists a symmetric matrix
$Z\in\mathbb{R}_+^{I\times I}$ such that $\Lloc(x) = x^T Z x$
is a local Lyapunov function with slack
$\varepsilon>0$. 
Then the $\varepsilon$-LRFS process is queue-length-stable if one of
the following conditions C1 and C2 is satisfied:
\begin{itemize}
	\item[C1.] The network is synchronized and
$Z_{i\ell}\neq 0$ only if buffers $i$ and $\ell$ are in the same component.
\item[C2.] 	
$Z_{i\ell}\neq 0$ only if buffers $i$ and $\ell$ are activity-interchangeable and
every buffer $i$ has an associated activity $j\in \J_i$ with $|\K_j|=1$.
\end{itemize}

Furthermore, if all routes
are bounded, then the LRFS process is queue-length-stable
if one of the conditions C1 or {C2$\,^\prime$} is satisfied, where
\begin{itemize}
\item[{C2$\,^\prime$}.] $Z_{i\ell}\neq 0$ only if buffers $i$ and $\ell$ are activity-interchangeable.
\end{itemize}
\end{theorem}
Theorem~\ref{thm:main} implies that if any maximal policy is
`quadratic' stable in a stochastic processing network without
routing under the external load $\brho$, then
$\varepsilon$-LRFS is stable in a stochastic processing network with
routing under nominal load $<\brho$ for some small
$\varepsilon>0$. Section~\ref{sec:overview} describes the main idea
of the proof, and a full proof is presented in Section
\ref{sec:pfthmmain}.

We remark that the requirement of activity-interchangeable buffers
in Condition {\em C2} can be relaxed slightly. Our proof
of Theorem \ref{thm:main} also works
when the following relaxed condition {\em C3} replaces {\em C2}.
\begin{itemize}
\item[{\em C3}.] $Z_{i\ell}\neq 0$ only if for every $j\in \J_i$, there exists
$j^{\prime}\in \J_\ell$ such that $\K_{j^{\prime}}\subseteq \K_j$, and
vice versa (i.e., for every $j\in \J_\ell$, there exists
$j^{\prime}\in \J_i$ such that $\K_{j^{\prime}}\subseteq \K_j$).
Also, every buffer $i$ has an associated activity $j\in \J_i$ with $|\K_j|=1$.
\end{itemize}
A corresponding condition {\em C3$\,^\prime$} can also replace {\em C2$\,^\prime$}, where {\em C3$\,^\prime$} does not require the second part of {\em C3}.

\section{Proof Ideas for Theorem
\ref{thm:main}}\label{sec:overview}

In this section, we describe the main idea in the proof of Theorem
\ref{thm:main}. We first present it in a very special network,
the Rybko-Stolyar network \cite{RybkoStolyar}, which
allows us to summarize the main idea of the proof at a high level.
We subsequently describe the challenges that have to be
overcome to establish our result in the general case,
and discuss the feasibility of an approach based on fluid models.

\subsection{Rybko-Stolyar Network} 
\begin{figure}[h]
\begin{center}
\includegraphics[width=11cm,angle=0]{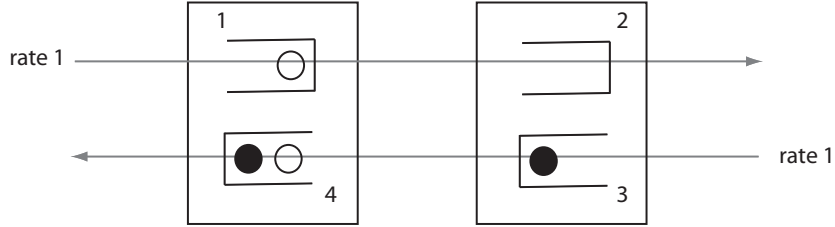}
\caption{Illustration of the Rybko-Stolyar Network. In this diagram,
there are four jobs in the system. Jobs are currently being
processed and waiting for service are colored black and white,
respectively.} \label{fig:RybkoStolyar}
\end{center}
\end{figure}

This network consists of four
activities associated to four different buffers and two processors.
The first processor is required for activities 1 and 4, and the
second for activities 2 and 3. Each activity decreases the remaining
service requirement of the job it is currently processing at unit
rate (i.e., $\beta_j=1$). Customers (or jobs) arrive at the first
and third buffers, and traverse the buffers deterministically in the
order $1\to2$ or $3\to4$. The service time is deterministic and equal
to $m_i$ for buffer $i$, and the external arrival processes (at the first and
third buffers) are independent Poisson processes with rate 1. The network is given in
Figure~\ref{fig:RybkoStolyar}, and a necessary condition for
stability is
\begin{equation}\label{eq:regionRybkoStolyar}
\rho_1=m_1+m_4<1\qquad \mbox{and}\qquad
\rho_2=m_2+m_3<1.\end{equation}

The network state can be described by
$X(t)=\{Q_i(t),V_i(t):i=1,2,3,4\}\in\Omega_X:=\mathbb{Z}_+^4\times\mathbb{R}_+^4$,
where $Q_i(t)$ is the queue length (i.e., the number of jobs waiting) of buffer $i$ at time $t$ and $V_i(t)$ is the
remaining service time of the job currently being processed in buffer
$i$ at time $t$. Hence $V_i(t)=0$ if no job is being processed in buffer $i$. We
define a norm through $|X(t)|=\sum_i (V_i(t)+Q_i(t))$. In this
network, it is known \cite{RybkoStolyar} that if processor 1
prioritizes buffer 4 and processor 2 prioritizes buffer 2, the
network can be unstable even when the necessary condition
\eqref{eq:regionRybkoStolyar} holds. On the other hand, LRFS
prioritizes buffer 1 and 3, so it reduces to the First Buffer First
Served policy. This is known to be stable under
the necessary stability condition \eqref{eq:regionRybkoStolyar}.
We next derive this stability result using our main idea.

Assuming LRFS and
\eqref{eq:regionRybkoStolyar}, we construct an appropriate Lyapunov
function $\Lglo:\Omega_X\to \mathbb{R}_+$ satisfying
\begin{eqnarray}\label{eq:driftRybkoStolyar}
E\left[\Lglo(X(t+1))-\Lglo(X(t))~\Bigg|~X(t)\right] &\leq&-\gamma
|X(t)|+\cC,\qquad\mbox{almost surely},\end{eqnarray}  where
$\gamma,\cC>0$ are some constants. By first taking expectations with
respect to the distribution of $X(t)$ and then integrating over $t$ on
both sides of the above inequality, we conclude that
\begin{eqnarray}
\limsup_{t\to\infty}\frac1t \int_{0}^{t}
E\left[|X(s)|\right]ds&\leq&
\limsup_{t\to\infty} \frac1{\gamma}\left(\cC - \frac{\int_t^{t+1} E[\Lglo(X(s))]ds-\int_0^1 E[\Lglo(X(s))]ds}t\right)\nonumber\\
&\leq&
\limsup_{t\to\infty} \frac1{\gamma}\left(\cC + \frac{\int_0^1 E\left[\Lglo(X(s))\right] ds}t\right)\nonumber\\
&=& \frac{\cC}{\gamma}~<~\infty.\label{eq:stabilityRybkoStolyar}
\end{eqnarray}

We now proceed toward constructing the `global' Lyapunov function
$\Lglo$ satisfying \eqref{eq:driftRybkoStolyar} based on a `local'
quadratic Lyapunov function $\Lloc$. To this end, we first discuss
how to construct the `local' quadratic Lyapunov function. Consider
a single-processor system with two buffers $a$ and $b$,
deterministic service times given by $m_a$ and $m_b$, and
independent Poisson arrival processes with rate $1$ at each buffer (i.e.,
the total rate is $2$). Hence, a necessary
condition for stability is
$$m_a+m_b<1.$$
Under a maximal (i.e., work-conserving)
scheduling policy, the workload $W(t)=\sum_i W_i(t)$
at time $t$ satisfies
$$W(t+1) ~=~ W(t)+\sum_{i=a,b} m_i\, A_i(t,t+1)-
\int^{t+1}_t {\bf 1}^+_{W(u)} du,$$
where
$W_i(t)$ is defined as the immediate workload at buffer $i$ as
in \eqref{eq:defW}, $A_i(s,t)$ is the number of jobs
arriving at buffer $i$ during the time interval $[s,t)$ so that
$E[A_i(s,t)]=t-s$ and we define
$${\bf 1}^+_x=\begin{cases} 1&\mbox{if}~x>0\\0&\mbox{otherwise}\end{cases}.$$
Hence, we have
\begin{equation}
E\left[W(t+1)-W(t)~\Big|~W(t)\right]~\leq~m_a+m_b-
{\bf 1}^+_{W(t)-1},\label{eq:singleserver0}
\end{equation}
where we use that ${\bf 1}^+_{W(u)}\geq {\bf 1}^+_{W(t)-1}$ for $u\in[t,t+1]$
since we assume deterministic service times. Using this, it
follows that for some finite constant $\cC$,
\begin{eqnarray}
E\left[ W(t+1)^2-W(t)^2~\Bigg|~W(t)\right]&=&
E\left[ 2W(t)(W(t+1)-W(t))~\Bigg|~W(t)\right]\notag\\
&~&\qquad+ E\left[
(W(t+1)-W(t))^2~\Bigg|~W(t)\right]\notag\\
&\leq&2W(t) E\left[W(t+1)-W(t)~\Bigg|~W(t)\right]+
\cC\notag\\&\leq&2W(t)\left(m_a+m_b-{\bf 1}^+_{W(t)-1} \right)+
\cC\notag\\
&<&-2(1-m_a-m_b)W(t)+\cC+2,\label{eq:singleserver}
\end{eqnarray}
where we use $E\left[A_i(t,t+1)^2\right]<\infty$ for the first
inequality. This shows that $W(t)^2=(W_1(t)+W_2(t))^2$ is a suitable
Lyapunov function for our `local' single-processor system under
the necessary stability condition $m_a+m_b<1$.

This observation on the single-processor system motivates the
following quadratic local Lyapunov function
$\Lloc:\mathbb{R}_+^I\to\mathbb{R}_+$ for the Rybko-Stolyar network:
\begin{equation}\label{eq:LlocRybkoStolyar}
\Lloc (x) ~=~ x^T Z x ~=~x^T \begin{pmatrix}
1&0&0&1\\0&1&1&0\\0&1&1&0\\1&0&0&1\end{pmatrix} x~=~(x_1+x_4)^2 +
(x_2 + x_3)^2,\qquad\mbox{for}~~x=(x_1,x_2,x_3,x_4).
\end{equation}
One can easily check that it satisfies \eqref{eq:localcondition}
with (small) slack $\varepsilon>0$ under the necessary stability
requirements $m_1+m_4<1$ and $m_2+m_3<1$. We propose the following
global Lyapunov function $\Lglo$:
$$
\Lglo(X(t))~=~\Lloc\big(W_1(t),V_2(t),W_3(t),V_4(t)\big)+ \xi
\Lloc\left(W_1(t),\widehat W_2(t),W_3(t),\widehat W_4(t)\right),
$$
where the new parameters $\xi$ and $\widehat W_i(t)$ shall be
defined explicitly. We remind the reader that our goal is to prove
\eqref{eq:driftRybkoStolyar}.

First, a similar calculation as for the single-processor case in
\eqref{eq:singleserver0} yields that under the LRFS policy,
\begin{eqnarray}
E\left[W_1(t+1)+V_4(t+1)-W_1(t)-V_4(t)~\Big|~X(t)\right]&\leq&
-{\bf 1}^+_{W_1(t)+V_4(t)-1}+m_1,\notag
\\E\left[W_3(t+1)+V_2(t+1)-W_3(t)-V_2(t)~\Big|~X(t)\right]&\leq&
-{\bf 1}^+_{W_3(t)+V_2(t)-1}+m_3.\notag
\end{eqnarray}
Hence, as for \eqref{eq:singleserver}, one can conclude that for
some constant $\cC$,
\begin{eqnarray}
E\left[\left(W_1(t+1)+V_4(t+1)\right)^2
-\left(W_1(t)+V_4(t)\right)^2~\Big|~X(t)\right]&\leq&
-2(1-m_1)(W_1(t)+V_4(t))+\cC,\notag\\
\label{eq1:globalRybkoStolyar}\\
E\left[\left(W_3(t+1)+V_2(t+1)\right)^2-\left(W_3(t)+V_2(t)\right)^2~\Big|~X(t)\right]&\leq&
-2(1-m_3)(W_3(t)+V_2(t))+\cC,\notag\\
\label{eq2:globalRybkoStolyar}
\end{eqnarray}
where the precise value of $\cC$ can be different from line to line.

We note that the sum
$\left(W_1(t)+V_4(t)\right)^2+\left(W_3(t)+V_2(t)\right)^2$ is not a
suitable choice for $\Lglo$ since it does not include $Q_2(t)$ and
$Q_4(t)$ (or $W_2(t)$ and $W_4(t)$). To address this issue, we
further use
$$\widehat W_2(t):=W_2(t)+m_2 Q_1(t)
\qquad\mbox{and}\qquad \widehat W_4(t):=W_4(t)+m_4Q_3(t).$$
We refer to $\widehat W_2$ and $\widehat W_4$ as the total workload
in buffer $2$ and $4$, respectively.
Using this notation, one finds that under the LRFS policy,
\begin{eqnarray}
E\left[W_1(t+1)+\widehat W_4(t+1)-W_1(t)-\widehat
W_4(t)~\Big|~X(t)\right]&\leq&
-{\bf 1}^+_{W_1(t)+W_4(t)-1}+m_1+m_4,\notag\\
E\left[W_3(t+1)+\widehat W_2(t+1)-W_3(t)-\widehat
W_2(t)~\Big|~X(t)\right]&\leq&
-{\bf 1}^+_{W_2(t)+W_3(t)-1}+m_2+m_3.\notag
\end{eqnarray}
The above equalities can be used to obtain `negative drift terms'
for $W_2(t)$ and $W_4(t)$, which are missing in
\eqref{eq1:globalRybkoStolyar} and
\eqref{eq2:globalRybkoStolyar}.
Namely, for some constant $\cC$, we obtain
\begin{eqnarray}
&&E\left[\left(W_1(t+1)+\widehat
W_4(t+1)\right)^2-\left(W_1(t)+\widehat
W_4(t)\right)^2~\Big|~X(t)\right]\notag\\
&&\qquad\qquad\leq~ 2\left(W_1(t)+\widehat W_4(t)\right)
E\left[W_1(t+1)+\widehat W_4(t+1)
-W_1(t)-\widehat W_4(t)~\Big|~X(t)\right]+\cC\notag\\
&&\qquad\qquad\leq~2\left(W_1(t)+\widehat
W_4(t)\right)\left(-{\bf 1}^+_{W_1(t)+W_4(t)-1}+m_1+m_4\right)+\cC\notag\\
&&\qquad\qquad\leq~-2(1-m_1-m_4)\left(W_1(t)+
W_4(t)\right)+2(m_1+m_4)m_4Q_3(t)+\cC+2\notag\\
&&\qquad\qquad\leq~-2(1-\rho_{\max})\left(W_1(t)+
W_4(t)\right)+2m^*\rho_{\max}
W_3(t)+\cC+2,\label{eq3:globalRybkoStolyar}
\end{eqnarray}
where $\rho_{\max}=\max\{\rho_1,\rho_2\}$ and
$m^*=\max_{i,j}\left\{\frac{m_i}{m_j}\right\}$. Similarly,
\begin{eqnarray}
&&E\left[\left(W_3(t+1)+\widehat
W_2(t+1)\right)^2-\left(W_3(t)+\widehat
W_2(t)\right)^2~\Big|~X(t)\right]\notag\\
&&\qquad\qquad\qquad\qquad\qquad\leq-2(1-\rho_{\max})\left(W_3(t)+
W_2(t)\right)+2m^*\rho_{\max}
W_1(t)+\cC+2.\label{eq4:globalRybkoStolyar}
\end{eqnarray}

Observe that there are positive terms $W_3(t)$ and $W_1(t)$ in
\eqref{eq3:globalRybkoStolyar} and \eqref{eq4:globalRybkoStolyar},
respectively.
The key idea behind our proof is that the positive terms can be canceled out by
appropriately summing \eqref{eq1:globalRybkoStolyar},
\eqref{eq2:globalRybkoStolyar}, \eqref{eq3:globalRybkoStolyar} and
\eqref{eq4:globalRybkoStolyar}. Indeed, we define the desired
Lyapunov function $\Lglo$ as
$$\Lglo(X)~=~\left(W_1+V_4\right)^2+\left(W_3+V_2\right)^2
+\xi \left(W_1+\widehat W_4\right)^2+\xi \left(W_3+\widehat
W_2\right)^2,$$ where we choose
$\xi=\frac{1-\rho_{\max}}{2m^*\rho_{\max}}$. Combining
\eqref{eq1:globalRybkoStolyar}, \eqref{eq2:globalRybkoStolyar},
\eqref{eq3:globalRybkoStolyar} and \eqref{eq4:globalRybkoStolyar},
we conclude that 
$$
E\left[\Lglo(X(t+1))-\Lglo(X(t))~\Bigg|~X(t)\right]~\leq~-2\xi(1-\rho_{\max})
|X(t)|+\cC.$$ This completes the proof of
\eqref{eq:driftRybkoStolyar}, and hence the desired stability
\eqref{eq:stabilityRybkoStolyar}.

\subsection{Beyond the Rybko-Stolyar Network}
The preceding subsection presents the main idea behind our
construction of a `global' Lyapunov function $\Lglo$ using a `local' Lyapunov function
$\Lloc$ (i.e., from the single-processor system) in the specific
example of the
Rybko-Stolyar network.
The construction of $\Lglo$ relies on summing $\Lloc$ terms inductively by exploring
certain maximality properties of the LRFS policy at each iteration.
In general networks there are several difficulties which do not arise in
the Rybko-Stolyar network, and
this section discusses the ideas and
arguments needed to overcome them.

A first challenge we have overcome
arises in networks with unbounded route lengths
(i.e., $\inf\{d\in\mathbb{N}:P^d=0\}=\infty$).
In that case, the above inductive procedure does not terminate. For this reason, we propose
a variant of the LRFS policy, the $\varepsilon$-LRFS policy,
which 
occasionally processes a job with the largest counter.
Intuitively speaking, this additional mechanism in
$\varepsilon$-LRFS can control the jobs
with large counters, whereas LRFS cannot.

A second challenge we have surmounted is that
the construction of $\Lglo$ in the
Rybko-Stolyar network starts from a simple local Lyapunov function
in a single-server system, but it is not clear
whether similar arguments go through for general local Lyapunov functions and stochastic
processing networks.
We require Condition {\em C2} to resolve this issue.
It is readily seen that the local Lyapunov function \eqref{eq:LlocRybkoStolyar} used
in the Rybko-Stolyar network satisfies this condition. The condition
can be relaxed under some additional conditions on the arrival processes
and service time distributions. For example, in synchronized networks,
Condition {{\em C1}} can be used instead of {\em C2}.

A further challenge in the general case relates to the definition of $\widehat W_i(t)$.
In the Rybko-Stolyar network, it is the sum of workloads along a path of buffers, with $i$ as the last
buffer. This definition only applies to networks with deterministic routing.
In the general case we use several notions of total workload.
To allow for stochastic routing, we construct a new process $\{Y(t)\}$ from $\{X(t)\}$
with deterministic routing. This process is essentially identical to $\{X(t)\}$,
but we enlarge the state space to incorporate routing information.
We construct a Lyapunov function $\Lglo$ for $\{Y(t)\}$, which we use to establish the stability of
$\{Y(t)\}$ and hence the stability of $\{X(t)\}$.

In summary, we construct the Lyapunov function
$\Lglo$ for general networks as the sum of three parts:
\begin{equation*}
	\Lglo(Y(t))~=~\sum_{c=1}^D \left(\frac{\upsilon}{2\mathcal C}\right)^c \Z\left(\widehat{\bW}_{\leq c}(t)\right)
+\frac{2\mathcal C}{\beta_{\min}}
\,\|\bV(t)\|_2^2+\frac{\xi}{2\cC}\,\mathcal G(Y(t))^2.\end{equation*}
The specific notation used here is not important;
we refer readers to Section \ref{sec:pfthmmain} for the definitions used.
To prove stability, we need to argue that this function satisfies a so-called negative-drift condition.
The first term, i.e., the finite sum,
comes from the inductive construction under LRFS,
appropriately truncated. For the Rybko-Stolyar network, this is the only part we need.
The first part produces the desired negative drift for jobs with low
counters, but it gives a positive drift in terms of remaining service requirements $\bV(t)$
as a by-product (albeit not in the Rybko-Stolyar network under the assumptions of the preceding
subsection).
The second term in our Lyapunov function ($\|\bV(t)\|^2_2$) has a negative drift
and compensates the positive drift incurred by the first term.
The third term in our Lyapunov function ($\mathcal G(Y(t))^2$) controls
the high-counter jobs under the
mechnism which is present in the $\varepsilon$-LRFS policy but not in the LRFS policy
(step 3 in Definition~\ref{def:epLRFS}).
This additional mechanism allows us to establish a negative drift for the last term.
By appropriately weighing each of the three terms, we derive the desired
negative drift condition for the Lyapunov function $\Lglo$.

\subsection{Connection with Fluid Models}
\label{sec:connectionfluid}

As mentioned in the previous subsection, our approach relies on an inductive argument
based on job counters. For the Rybko-Stolyar network and more generally for multiclass networks,
fluid models can be used to give relatively simple proofs of our results.
Thus, a more detailed discussion on the connection with fluid models together with its pros and cons is warranted.

The fluid approach consists of two main steps. In the first step, by scaling time and space, one proves convergence
of the queueing process to the solution of a system of deterministic equations known as the {\em fluid model}.
In the second step, one proves that this fluid model is stable, i.e., that it eventually reaches the origin.
Stability of the fluid model can be established through the construction of a Lyapunov function for the fluid model,
or in some cases one can obtain fluid stability through direct methods such as induction.
Once fluid stability has been established,
one can apply general theorems to deduce that the stochastic model is also stable (in a certain sense), see
for instance \cite{Dai95}.

It might be possible to establish existence of a fluid model and to prove that the stochastic model converges to the fluid model
in the setting of the present paper, and
it can be expected that our `global' Lyapunov function should work to prove fluid stability.
Comparing our approach with this proof strategy, a disadvantage of the fluid model is that
one needs to establish convergence to the fluid model, while a disadvantage of our approach is that
we have to keep track of detailed state information such as residual service times.

Another possible approach to establish fluid stability is to use an inductive argument, which may
seem particularly attractive given our construction of job counters
and the suitability of a induction argument in existing work on fluid models \cite{Dai96}.
However, this approach has inherent challenges.
The base step in an inductive approach could use the `local' Lyapunov function
$\Lloc$ to argue that the fluid level of
jobs with counter 1 vanishes after some finite time $T_1$.
It would then use $\Lloc$
to argue that the fluid level of jobs with counter 2 vanishes after some finite time $T_2>T_1$, and so forth.
To carry out this argument, one has to show that $\Lloc$
satisfies a certain negative-drift condition under the assumption that high-priority counter 1 jobs
vanish on a fluid scale.
The latter only yields a guarantee on the `average' or `long-run' behavior
of the jobs with counter 1, whereas one needs
`short-term' network state information to establish the negative-drift condition for jobs with counter 2.
Indeed, under our scheduling policy, jobs with counter 1 (even when vanishing on a fluid scale)
can significantly influence the dynamics of jobs with counter 2 depending on the complexity of the network. 
Therefore, the base of the induction approach is too weak to be used in the induction step for general networks since
one needs more detailed information than the time-average given by the fluid approach.

In special cases such as multiclass networks, one may not need quadratic Lyapunov functions and
it may be possible to establish the stability of our counter-based policy using fluid induction
without quadratic Lyapunov functions. However, in general (e.g., for networks of switches), we
need quadratic Lyapunov functions since they are the only
available tool to establish stability for single-hop networks.

\section{Examples}\label{sec:application}
In this section, we provide applications of Theorem \ref{thm:main}
to various special stochastic processing networks. We consider
parallel server networks (including multiclass queueing
networks) in Section \ref{sec:mcqn} and communication networks (including
wireless networks and networks of input-queued switches) in Section
\ref{sec:swn}. They are examples of non-synchronized and synchronized
networks, respectively. In all of these important examples, suitable local
Lyapunov functions are easy to find. 

\subsection{Open Multiclass Queueing Networks and Parallel Server
Networks}\label{sec:mcqn}

In this section, we consider special
stochastic processing networks
known as {\em parallel server networks}.
These networks are characterized by the following assumption.

\begin{figure}[h]
\begin{center}
\includegraphics[width=15cm,angle=0]{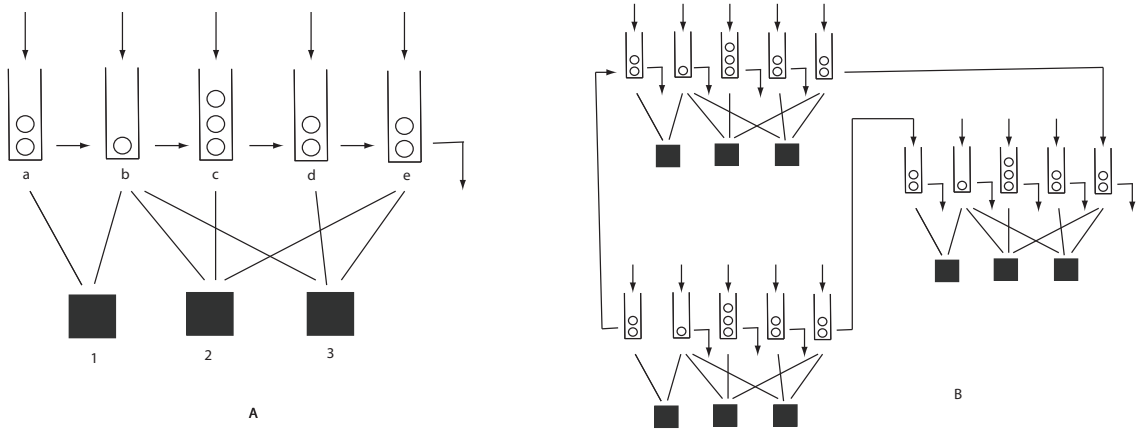}
\caption{Two examples of parallel server networks. They do not satisfy
Assumption {\em A2}. The leftmost diagram
illustrates a parallel server network with five
buffers $\{a,b,c,d,e\}$, three processors $\{1,2,3\}$ and eight
activities $\{(a,1),(b,1),(b,2),(b,3),(c,2),(d,3),(e,2),(e,3)\}$.
Once a job in buffer $a,b,c,d$ completes its service requirement, it
joins buffer $b,c,d,e$ (i.e., $P_{ab}=P_{bc}=P_{cd}=P_{de}=1$),
respectively. Once a job in buffer $e$ completes its service
requirement, it leaves the network. The rightmost diagram illustrates a
parallel server network consisting of three `local' parallel server systems
where jobs are routed between local systems.
} \label{fig:psn}
\end{center}
\end{figure}

\begin{itemize}
\item[\em A1.] Each activity is processed by exactly
one processor and processes exactly one buffer, i.e.,
$$|\K_j|=1,\qquad\mbox{ for all}~~j\in \J.$$
\end{itemize}
Figure \ref{fig:psn} illustrates the relations between buffers,
activities and processors in parallel server networks.
Our notion of `parallel server network' generalizes the well-studied parallel server
systems \cite{harrison0} by adding stochastic routing dynamics
between buffers. It also includes open multiclass
queueing networks \cite{Dai95} as a special case, which additionally
require
\begin{equation}
|\J_i|=1,\qquad\mbox{ for all}~~i\in
\I.\label{eq:conmqn}\end{equation} In open multiclass queueing
networks, buffers and activities are in one-to-one correspondence
and they are referred to as {\em classes}. The Rybko-Stolyar
in Section \ref{sec:overview} is an instance of open multiclass
queueing networks.

\begin{figure}[h]
\begin{center}
\includegraphics[width=14cm,angle=0]{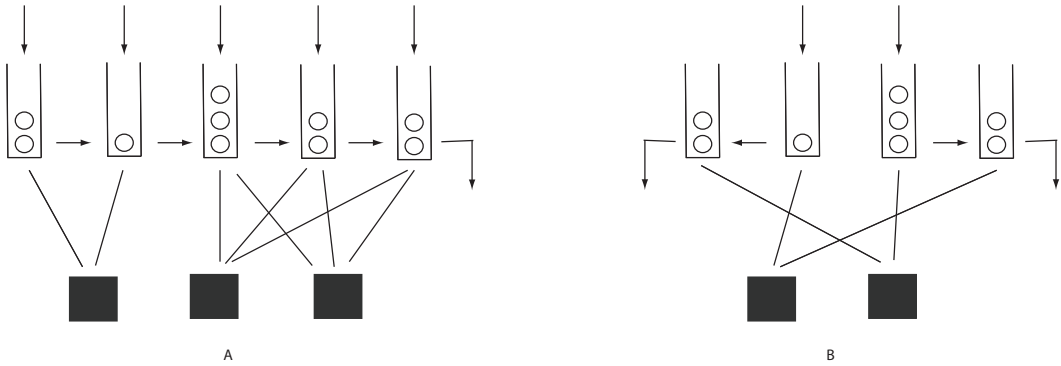}
\end{center}
\caption{Two examples of parallel server networks satisfying
Assumption {\em A2}. The rightmost diagram illustrates the
Rybko-Stolyar network described in Section \ref{sec:overview}.}
\label{fig:psn2}
\end{figure}

A parallel server network
naturally defines a bipartite graph $(\I,\K,\J)$ such that each
activity in $\J$ defines an edge between buffers $\I$ and processors
$\K$. Requirement \eqref{eq:conmqn} of open multiclass queueing
networks imposes the additional restriction that each
vertex in $\I$ has degree one. We further consider the following strengthening of
Assumption {\em A1}.
\begin{itemize}
\item[{\em A2}.] $(\I,\K,\J)$ is a
union of disjoint complete bipartite graphs
$\{(\I^{(h)},\K^{(h)},\J^{(h)}):h\in\cH\}$, i.e.,
$$\cup_{h\in \cH} \I^{(h)}=\I,\qquad
\cup_{h\in\cH} \K^{(h)}=\K,\qquad \mbox{and}\qquad
\cup_{h\in\cH} \J^{(h)}=\J.$$
\end{itemize}
This assumption implies that
two buffers in the same component are 
activity-interchangeable, even though
they may differ with respect to routing, external arrivals or service requirements.
One can easily check that open multiclass queueing networks always
satisfy this assumption, while the parallel server networks
in Figure \ref{fig:psn} do not. Assumption {\em A2} is useful
because it enables us to establish a necessary
condition for stability and it allows us to find a suitable local Lyapunov
function satisfying Condition {\em C2} of Theorem \ref{thm:main}.
However, Theorem \ref{thm:main} is applicable to
general networks as long as one can find a `good' local Lyapunov
function satisfying Condition {\em C2}. Figure \ref{fig:psn2}
gives examples of parallel server networks satisfying Assumption {\em A2}.

\paragraph{Necessary condition for stability.} We now aim to obtain a necessary
condition for stability of a parallel server network.
Under Assumption {\em A2}, a necessary condition to stabilize the
network is that for every $n$,
\begin{equation}\label{eq:stableregionpsn}
\rho^{(h)}:=\sum_{i\in\I^{(h)}} \rho_i~<~\beta^{(h)}:=\sum_{j\in\J^{(h)}}
\beta_j.
\end{equation}
It is clear that the above condition is required for stability since
$\rho^{(h)}$ and $ \beta^{(h)}$ describe the total nominal load
and the maximum processing rate, respectively, at the local
component $(\I^{(h)},\K^{(h)},\J^{(h)})$.

\paragraph{Local Lyapunov function.}
As in Section~\ref{sec:overview}, the single-processor example is
the main building block. We define the local Lyapunov function as
\begin{equation}\label{eq:mulloc} \Lloc(x) = \sum_{h\in\cH}
\left(\sum_{i\in \I^{(h)}} x_{i}
\right)^2\qquad\mbox{for}~~x=(x_1,\dots,x_I).
\end{equation} We now show that this function satisfies
(\ref{eq:localcondition}) with some slack $\varepsilon\geq 0$ as
long as the necessary
condition \eqref{eq:stableregionpsn} for stability is satisfied. 
For given vectors $\bw=[w_i]\in\mathbb{R}_+^I$ and
$\bu=[u_j]\in\mathcal M(\bw)$,
maximality implies that, on writing  
$w^{(h)}=\sum_{i\in \I^{(h)}} w_{i}$,
$$ s^{(h)}=s^{(h)}(\bu)
~:=~\sum_{j\in\J^{(h)}}u_{j}\beta_j~\geq~ \sum_{j\in\J^{(h)}}
{\bf 1}^+_{
w^{(h)}}\beta_j~=~{\bf 1}^+_{ w^{(h)}} \beta^{(h)},$$ where
we use Assumption {\em A2} and we
recall that ${\bf 1}^+_x=1$ if $x>0$, and ${\bf 1}^+_x=0$ otherwise.
Thus, we have
\begin{eqnarray*}
\Lloc(\bw+\brho+\varepsilon\bphi-s(\bu))-\Lloc(\bw)&=& \sum_{h\in\cH}
\left[\left( w^{(h)}+ \rho^{(h)}+\varepsilon\, m^{(h)}-
s^{(h)}\right)^2-\left(
w^{(h)}\right)^2\right]\\
&\leq&\cC+2\sum_{h\in\cH}
w^{(h)}\left(\rho^{(h)}+\varepsilon\, m^{(h)}- s^{(h)}\right)\\
&\leq&\cC+2\sum_{h\in\cH}
w^{(h)}\left(\rho^{(h)}+\varepsilon\, m^{(h)}-
{\bf 1}^+_{w^{(h)}}\beta^{(h)}\right)\\
&=&\cC+2\sum_{h\in\cH}
w^{(h)}\left(\rho^{(h)}+\varepsilon\, m^{(h)}-\beta^{(h)}\right),\end{eqnarray*}
where $\cC$ is some constant and we define
$m^{(h)}:=\sum_{i\in\I^{(h)}}m_i.$ Therefore, $\Lloc$ is a local
Lyapunov function with slack $\varepsilon$ for 
$$0\leq \varepsilon<\min_{h\in\cH} \frac{ \beta^{(h)}-\rho^{(h)}}{m^{(h)}},$$
where the right-hand side is positive if
\eqref{eq:stableregionpsn} holds. 

\paragraph{Stability of LRFS policies.}
We now formulate the main results of this paper for open multiclass
networks and parallel server networks. Under Assumption {\em A2}, the local
Lyapunov function \eqref{eq:mulloc} satisfies Condition {\em C2} of
Theorem \ref{thm:main}.
Therefore, we obtain the
following proposition as a corollary. We remind the reader that open multiclass
queueing networks are special instances of parallel
server networks, and that Assumption {\em A2} automatically holds for these
networks.

\begin{proposition}\label{cor:psn}
If a stochastic processing network satisfies Assumption {\em A2}
with $ \rho^{(h)}< \beta^{(h)}$ for all $n$, then
\begin{itemize}
\item The $\varepsilon$-LRFS process is queue-length-stable for any
$\varepsilon \in\left(0,\min\limits_{h\in\cH}
\frac{\beta^{(h)}-\rho^{(h)}}{m^{(h)}}\right).$
\item The LRFS process is queue-length-stable if all routes are bounded in length.
\end{itemize}
\end{proposition}
We note that the $\varepsilon$-LRFS policy admits a simpler description
in a stochastic processing network satisfying Assumption {\em A2},
since a job can be processed by any processor in the partition, i.e.,
$\Sigma$ in Definition \ref{def:epLRFS} is non-empty whenever
a processor is idle and capable of processing a job.
Indeed, the $\varepsilon$-LRFS policy reduces to
the following work-conserving randomized priority policy:
whenever a processor $k$ is idle at time $t$ and there are jobs capable of being processed by $k$,
\begin{itemize}
\item Process a job with
the smallest counter with probability $1-\varepsilon \left(1-{\bf 1}^+_{\mathcal T^{(h)}(t)}\right)$,
otherwise process a job with the largest counter.
\item Set $\mathcal T^{(h)}(t)=1$ if $\mathcal T^{(h)}(t)=0$,\end{itemize}
where $\I^{(h)}$ is the (local) component of buffers associated with
processor $k$.
Proposition \ref{cor:psn} implies that the $\varepsilon$-LRFS policy
can achieve `almost' the full capacity region \eqref{eq:stableregionpsn}
by choosing a small $\varepsilon>0$.

Our proof of Proposition~\ref{cor:psn} provides a different proof for some results that have been established
using fluid model techniques.
For example, in reentrant lines, the $\varepsilon$-LRFS policy for $\varepsilon=0$ is
identical to the well-known
First Buffer First Served (FBFS) policy. Our proposition implies that the FBFS policy is throughput optimal
in all reentrant lines, which has been proved
originally in \cite{Dai96}.

\subsection{Communication Networks}\label{sec:swn}

We now consider examples of synchronized stochastic processing
networks described in Section \ref{sec:defSPN}, i.e.,
$m_i=\beta_j=1$, for all $i\in \I, j\in \J$.
In particular, we
consider the following additional assumption on synchronized stochastic
processing networks.
\begin{itemize}
\item[{\em B1.}] Each buffer has exactly one associated activity, i.e.,
$$\J_i=\{j_i\}\quad \mbox{for all}~~ i\in \I.$$
Hence, we write $\K_i:=\K_{j_i}$.
\end{itemize}
We again remark that Assumption {\em B1} facilitates a
suitable local Lyapunov function for Theorem \ref{thm:main}.
However, even if Assumption {\em B1} does not hold, Theorem \ref{thm:main} is
applicable to synchronized stochastic processing networks as long as one can find 
a `good' local Lyapunov function.
Synchronized stochastic processing networks satisfying Assumption {\em B1}
include various communication network models of unit-sized packets: 
networks of input-queued
switches \cite{McKeown,Dai00}, wireless network models with
primary interference constraints \cite{sanghavi2007distributed}
and independent-set interference constraints \cite{jinwoo}.
We refer the corresponding references for detailed
descriptions of the network models. As a concrete example,
we write out the
details of the wireless network model with primary interference
constraints.

\begin{figure}[h]
\begin{center}
\includegraphics[width=15cm,angle=0]{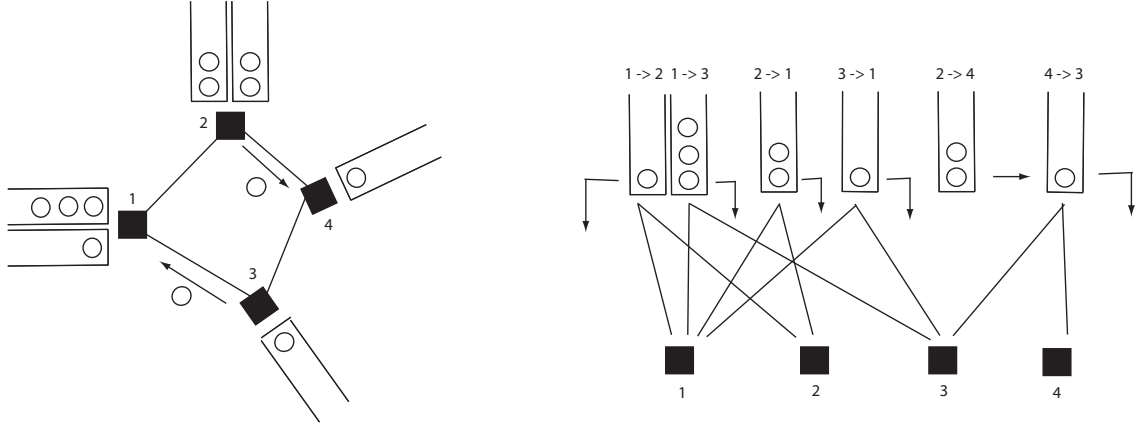}
\caption{Example of a wireless network with primary interference
constraints with four nodes $\{1,2,3,4\}$ and five paths $\{1\to
2,2\to1,  1\to 3,3\to 1,2\to 4\to 3\}$.
There is a buffer on each (directed) edge on each path, i.e., six buffers in total.
Unit-size packets arrive at the ingress buffer (i.e., the first node) of
each path. In the leftmost diagram,
two links $2\to 4$ and $3\to 1$ are
transmitting (unit-sized) packets.
Once a packet is transmitted, it
leaves the network if it arrives at the destination node (i.e., the last node on
its path). The rightmost diagram
illustrates the corresponding stochastic processing network,
where the relation
between buffers, activities and processors induces a hypergraph
(each buffer requires two processors)}.
 \label{fig:lsn}
\end{center}
\end{figure}

\paragraph{Wireless networks with primary interference
constraints.} Consider a network of $n$ nodes represented by
$V=\{1,\dots,n\}$ and a set of directed paths $\{P_1,P_2,\dots\}$.
Unit-size packets arrive at the ingress node of each path as
per an exogenous arrival process. Assume that the network is
synchronized, i.e., each packet departs from a node at time $t\in
\mathbb Z_+$ and arrives at the next node on its route at time
$t+1$. The primary interference constraint means that each node can
either send or receive (it cannot do both) one packet at the time. 
A scheduling policy (or algorithm) decides which packets transmit at
each (discrete) time instance. Figure \ref{fig:lsn}
illustrates a wireless network of four nodes with primary interference
constraints.

\paragraph{Necessary condition for stability.}
As in the parallel server networks, one can obtain the following necessary condition to
stabilize a stochastic processing network satisfying Assumption {\em B1}: for all $k\in \K$,
\begin{equation}\label{eq:stableregionsn}
\rho^{(k)}:=\sum_{i\in\I:k\in \K_i} \rho_{i}~<~1.\end{equation} This
is because Assumption {\em B1} implies that each buffer has at most
one associated activity, and the processing rate is $1$.

\paragraph{Local Lyapunov function.} We consider the following
local Lyapunov function:
$$\Lloc( x)=\sum_{i\in\I} \sum_{k\in \K_i}\sum_{\ell\in\I:k\in \K_\ell} x_i x_\ell,
\qquad\mbox{for}~~ x=( x_1,\dots, x_I).$$ We now proceed toward proving that
condition (\ref{eq:localcondition}) holds. 
For given vectors $\bw=[w_i]\in\mathbb{R}_+^I$
and $\bu=[u_j]\in\mathcal M(\bw)$,
maximality implies that 
$$\sum_{k\in \K_i}\sum_{\ell\in\I:k\in \K_\ell}u_\ell~\geq~
{\bf 1}^+_{w_i}.$$
Thus, we have that, on writing $\nu_i=\rho_i+\varepsilon m_i$,
\begin{eqnarray*}
\Lloc(\bw+\brho+\varepsilon\bphi-s(\bu))-\Lloc(\bw)&\leq&\cC+2
\sum_{i\in \I}w_i \left(\sum_{k\in \K_i}\sum_{\ell\in\I:k\in
\K_\ell}\nu_\ell-
\sum_{k\in \K_i}\sum_{\ell\in\I:k\in \K_\ell}u_\ell\right)\\
&\leq&\cC+2 \sum_{i\in \I}w_i \left(\sum_{k\in
\K_i}\sum_{\ell\in\I:k\in
\K_\ell}\nu_\ell-{\bf 1}^+_{w_i}\right)\\
&=&\cC+2 \sum_{i\in \I}w_i \left(\sum_{k\in \K_i}\sum_{\ell\in\I:k\in
\K_\ell}\nu_\ell-1\right),\end{eqnarray*} where $\cC$ is some
constant. Therefore, (\ref{eq:localcondition}) holds if
$\nu^{(k)}:=\sum_{i\in\I:k\in \K_i} \nu_{i}<\frac1{\max_{i\in\I}
|\K_i|}$ for all $k\in \K$. This is equivalent to
$$0\leq \varepsilon<\min_{k\in\K}\frac{\frac1{\max_{i\in\I}
|\K_i|}-\rho^{(k)}}{m^{(k)}},$$ where we define
$m^{(k)}:=\sum_{i\in\I:k\in \K_i} m_{i}$. The above interval is
non-empty as long as $\rho^{(k)}<\frac1{\max_{i\in\I} |\K_i|}$ for all
$k\in\K$.

\paragraph{Stability of LRFS policies.}
We now state the main result for synchronized stochastic networks
satisfying Assumption {\em B1}.
Since the network is synchronized, we obtain the following proposition as a corollary
of Theorem \ref{thm:main}.
\begin{proposition}\label{cor:sn}
If a stochastic processing network is synchronized and satisfies
Assumption B1 with $\rho^{(k)}< \frac1{\max_{i\in\I} |\K_i|}$
for all $k\in\K$, then
\begin{itemize}
\item The $\varepsilon$-LRFS process is queue-length-stable for any
$\varepsilon \in\left(0,\min_{k\in\K}\frac{\frac1{\max_{i\in\I}
|\K_i|}-\rho^{(k)}}{m^{(k)}}\right).$
\item The LRFS process is queue-length-stable if all routes are bounded in length.
\end{itemize}
\end{proposition}
Proposition \ref{cor:sn} implies that the $\varepsilon$-LRFS policy can
achieve a $\frac1{\max_{i\in\I} |\K_i|}$ fraction of the capacity
region \eqref{eq:stableregionsn}. For networks of input-queued
switches \cite{McKeown,Dai00} and wireless networks with
primary interference constraints \cite{sanghavi2007distributed}, it
is easy to see that $\max_{i\in\I} |\K_i|=2$, and hence the
$\varepsilon$-LRFS policy achieves 50\% of the capacity region. For
wireless networks with general independent set constraints
\cite{jinwoo}, $\max_{i\in\I} |\K_i|$ is the maximum number of
interfering neighbors, i.e., the maximum degree of the underlying
interference graph.

In large-scale networks, distributed scheduling schemes 
with low complexity have gained much attention recently, even though
they usually perform worse than centralized ones with high complexity.
We remark that the 50\% throughput result of a greedy scheduling algorithm
has previously been established for input-queued switches
(i.e., no routing between buffers) by Dai et al.\ \cite{Dai00}.
Proposition \ref{cor:sn} generalizes this to `networks' of input-queued switches
operated under stochastic routing between buffers (or local switches).
In wireless networks with primary interference constraints,
Wu et al.\ \cite{wu2007scheduling} establish
the 50\% throughput result of the LRFS policy
assuming deterministic routing (i.e., fixed routes between nodes),
while Proposition \ref{cor:sn} allows 
stochastic routing.

\section{Proof of Theorem \ref{thm:main}}\label{sec:pfthmmain}

To prove the desired stability of the $\varepsilon$-LRFS process
$\{X(t)\}$, we construct a new process $\{Y(t)\}$, which is almost identical
to $\{X(t)\}$, but it has a larger state space.
The main idea of the proof is to construct
a Lyapunov function for the `larger' process $\{Y(t)\}$, which implies
the stability of $\{Y(t)\}$ and therefore the stability of $\{X(t)\}$.

The description of $\{Y(t)\}$ is as
follows. Consider the stochastic processing network
setup in Section \ref{sec:defSPN}, and let
$\cP=\{P_1,P_2,\dots\}$ be the collection of all possible paths of
buffers (allowing repetitions) of length at most $D$ (i.e.,
$|\cP|\leq (I+1)^D$), where $D\in \mathbb{N}$ is some finite
constant to be determined later. 
We assume that when a job enters the network, it pre-determines the
first $D$ buffers on its route. After being processed from these $D$
buffers, jobs perform the usual stochastic routing as described in
Section \ref{sec:defSPN}. Let $i_{nm}\in\I$, $Q_{nm}(t)\in
\mathbb{Z}_+$ and $V_{nm}^j(t)\in \mathbb{R}_+$ denote the $m$-th
buffer on path $P_n$, the number of jobs waiting in buffer $i_{nm}$
and the remaining service requirement of the job in buffer $i_{nm}$
being processed by activity $j$ at time $t$, respectively. If
activity $j$ is not processing a job in buffer $i_{nm}$, then
$V_{nm}^j(t)=0$. Furthermore, as before, let $Q_{i,c}(t)$ be the queue length
(i.e., the number of jobs waiting for service, excluding those being processed) with counter $c$ in buffer $i$ at time $t$.
$V_{i,c}^j(t)$ is defined to be the remaining service requirement of the job with counter $c$ in
buffer $i$ being processed by activity $j$ at time $t$.
We then define
\begin{eqnarray*}
	Y(t)&=&\left[Q_{nm}(t), Q_{i,c}(t), V_{nm}^j(t), V_{i,c}^j(t),
\mathcal T^{(h)}(t)\right]
\in\Omega_Y:=\mathbb{Z}_+^{\infty}\times
\mathbb{R}_+^{\infty}\times [0,1]^H,\end{eqnarray*}
where $n,m,c,h$ are positive integers such that $n\leq (I+1)^D$, $m\leq D$, $c> D$ and $h\leq H$.

As for $\{X(t)\}$, we impose the convention that $\{Y(t)\}$ has right-continuous sample paths.
We define the norm $|Y(t)|$ through
$$|Y(t)|~=~\sum_{n\leq (I+1)^D,m\leq D} Q_{nm}(t)+\sum_{i\in\I, c>D} Q_{i,c}(t)
+\sum_{n\leq (I+1)^D,m\leq D,j\in \J} V^j_{nm}(t)
+\sum_{i\in\I, c>D,j\in\J} V^j_{i,c}(t).
$$

One can define a natural projection $P$
such that the distribution of $P(Y(t))$ is identical to that of
$X(t)$ and $|P(Y(t))|=|X(t)|$ given that $Y(0)$ is drawn appropriately from
the preimage $P^{-1}(X(0))$ of $X(0)$. Intuitively speaking, $\{Y(t)\}$ tosses
random coins to determine routes in advance, and since the scheduling
decisions of $\varepsilon$-LRFS are independent of these coins, the
natural projection of $\{Y(t)\}$ ignoring these pre-determined coin flips
provides exactly the dynamics of $\{X(t)\}$. Hence, it suffices to prove
that the (bigger) process $\{Y(t):t\in\mathbb{R}_+\}$ is
queue-length-stable, i.e.,
\begin{equation}\label{eq:boundY}
\limsup_{t\to\infty}\frac1t\int_0^t E\left[|Y(s)|\right]ds<\infty,
\end{equation}
for any given initial state $Y(0)\in \Omega_Y^*:=\{Y\in \Omega_Y:|Y|<\infty\}$. In
essence, this follows from the following proposition, which is proved in
Section \ref{sec:pflem2}.

\begin{proposition}\label{lem2}
If the conditions of Theorem \ref{thm:main} hold, then there exist 
constants $N,T,D,\zeta,\cC\in(0,\infty)$ and a Lyapunov function
$\mathcal L_{\text{\em
global}}:\Omega_Y\to\mathbb{R}_+\cup\{\infty\}$ such that for all
$t\geq N$, 
\begin{eqnarray*}
E\left[\mathcal L_{\text{\em
global}}(Y(t+T))~|~\cF(t)\right]&\leq&\mathcal L_{\text{\em
global}}(Y(t))-\zeta |
Y(t)|
+\mathcal C,\qquad\mbox{almost surely},
\end{eqnarray*}
and $\sup_{Y\in\Omega_Y^*}{\mathcal L_{\text{\em
global}}(Y)}/{\mathcal |Y|^2}<\infty$, where we define
the filtration $\{\cF(t):t\ge 0\}$ by
$$\mathcal F(t):=\sigma\{Y(s):0\leq s\leq t\}.$$
\end{proposition}

Now we describe how Proposition \ref{lem2} implies \eqref{eq:boundY},
and hence the conclusion of Theorem \ref{thm:main}.
First one can observe that $E\left[|Y(t)|^2\right]<\infty$ since
$Y(0)\in\Omega_Y^*$ and we assume bounded second moments on
arrivals and service times. Since $\sup_{Y\in\Omega_Y^*}{\Lglo(Y)}/{\mathcal |Y|^2}<\infty$ from Proposition \ref{lem2}, it follows that
\begin{equation}\label{eq:boundsec}
E\left[\Lglo(Y(t+T))\right],E\left[\Lglo(Y(t))\right],E\left[|
Y(t)|\right]<\infty.\end{equation}
Combining Proposition \ref{lem2} and \eqref{eq:boundsec} yields that for all
$t\geq N$,
$$E\left[\Lglo(Y(t+T))\right]~\leq~E\left[\Lglo(Y(t))\right]-\zeta
E\left[| Y(t)|\right]+\mathcal C.$$

Therefore, we have that for $t\geq N$,
\begin{eqnarray*}
\int^{T+t}_{T+N} E\left[\Lglo(Y(s))\right]ds&=&\int^t_N
E\left[\Lglo(Y(s+T))\right]ds\\
&\leq&\int^t_N E\left[\Lglo(Y(s))\right]ds-\zeta \int^t_N E\left[|
Y(s)|\right]ds+\mathcal C (t-N),\end{eqnarray*} which implies that
for $t\geq T+N$,
$$\frac{\zeta}{t-N}\int^t_N E\left[| Y(s)|\right]ds~\leq~\mathcal C+\frac1{t-N}\int^{T+N}_N
E\left[\Lglo(Y(s))\right]ds.$$ The right-hand side of the
above inequality converges to $\cC$ as $t\to\infty$. This leads to
the desired conclusion \eqref{eq:boundY}.

\subsection{Proof of Proposition \ref{lem2}}\label{sec:pflem2}

The choice of $N$ in Proposition \ref{lem2} comes from our assumption
on the external arrival processes in Section \ref{sec:defSPN},
namely, that there exists some $N<\infty$ such that for all $i\in \I$ and
$t\geq N$,
\begin{equation*}
E\Big[A_i(t,t+1)~\big|~\{A_1(0,s),\dots,A_I(0,s):0\leq s\leq
t\}\Big]~\leq~\alpha_i.\end{equation*} For notational convenience, we
assume $N=0$ in the proof of Proposition \ref{lem2}.
Namely, we assume that for all $t\geq 0$, 
\begin{equation}\label{eq:arrivalassum}
E\Big[A_i(t,t+1)~\big|~\cF(t)\Big]~\leq~\alpha_i<\infty.\end{equation}
All the proof arguments are applicable to the general
case $N>0$.

We first define some further notation in Section
\ref{sec:not}. The skeleton of the proof of Proposition \ref{lem2}
is described in Section \ref{sec:skeleton}, and it uses three key lemmas.
The proofs of these lemmas follow.

\subsubsection{Notation}\label{sec:not}
The quantities defined below are simple functions of the network
state $Y(t)$. 
We use bold symbols to denote vectors of quantities, e.g.,
$\bQ_{c}(t)=[Q_{i,c}(t)]$ and $\bW_{\leq c}(t)=[W_{i,\leq c}(t)]$.

\paragraph{Queues.}
For $c\in \mathbb{N}$, let $Q_{i,c}(t)$, $Q_{i,<c}(t)$ and
$Q_{i,\leq c}(t)$ be the number of waiting jobs with counter $c$,
$<c$ and $\leq c$ in buffer $i$ at time $t$, respectively. That is,
we set
\begin{eqnarray*}
Q_{i,c}(t)~=~\sum_{n=1}^{(I+1)^D}Q_{nc}(t)\,
I_{\{i\}}(i_{nc})\qquad\mbox{for}~~c\leq D,\\
Q_{i,<c}(t)~=~\sum_{c^{\prime}=1}^{c-1}Q_{i,c^{\prime}}(t)\qquad
Q_{i,\leq c}(t)~=~\sum_{c^{\prime}=1}^{c}Q_{i,c^{\prime}}(t),
\end{eqnarray*}
where we recall that
$I_{\{i\}}(i_{nc})=1$ if $i_{nc}=i$ and $I_{\{i\}}(i_{nc})=0$ otherwise.
Furthermore, for $c\leq D$, we let $\widehat{Q}_{i,\leq c}(t)$ be
the number of waiting jobs with counter $\leq c$ `in or destined
for' buffer $i$ at time $t$. Namely,
\begin{eqnarray*}
\widehat{Q}_{i,\leq c}(t)&=&{Q}_{i,\leq
c}(t)+\sum_{n=1}^{(I+1)^D}\sum_{m=1}^c
\sum_{r=m+1}^cQ_{nm}(t)\,I_{\{i\}}(i_{nr}),
\end{eqnarray*}
where the first term `${Q}_{i,\leq
c}(t)$' and the second term `$\sum_{n=1}^{(I+1)^D}\cdots$'
on the right-hand side count the numbers of
waiting jobs currently in buffer $i$
and destined for buffer $i$, respectively.
\paragraph{Workloads.}
Let $V^j(t)$ denote the remaining service requirement
of the job being processed by activity $j$ at time
$t$, and $V_i(t)$ the total remaining service
requirement of the jobs being processed in buffer $i$ at time $t$
(multiple jobs can be processed from the same buffer by different
processors). Similarly, $V_{i,c}(t),V_{i,<c}(t),V_{i,\leq c}(t)$
stands for the total remaining service requirement of jobs with
counter $c$, $<c$, $\leq c$ being processed in buffer $i$ at time
$t$, respectively. We furthermore define the following quantities:
\begin{itemize}
\item[] if the network is not synchronized,
\begin{eqnarray*}
{W}_{i,<c}(t)&=&m_i{Q}_{i,<c}(t)+V_i(t)\\
{W}_{i,\leq c}(t)&=&m_i{Q}_{i,\leq c}(t)+V_i(t)\\
\widehat{W}_{i,\leq c}(t)&=&m_i\widehat{Q}_{i,\leq c}(t)+V_i(t).
\end{eqnarray*}
\item[] if the network is synchronized,
\begin{eqnarray*}
{W}_{i,<c}(t)&=&m_i{Q}_{i,<c}(t)+V_{i,<c}(t)\\
{W}_{i,\leq c}(t)&=&m_i{Q}_{i,\leq c}(t)+V_{i,\leq c}(t)\\
\widehat{W}_{i,\leq c}(t)&=&m_i\widehat{Q}_{i,\leq c}(t)+V_{i,\leq
c}(t).
\end{eqnarray*}
\end{itemize}
We use different definitions for these variables depending on whether the network
is synchronized or not, since the strategy of the proof differs in each case,
e.g., see Section \ref{sec:pflemkey1}.
We also note that in a synchronized network, $V_{i, c}(t),V_{i,<
c}(t),V_{i,\leq c}(t)\in\{0,1\}$ for $t\in\mathbb Z_+$.
The quantities $W_{i,<c}(t)$ and $W_{i,\le c}(t)$ are
(expected) immediate workloads, since they only involve work that is in buffer $i$
at time $t$. The quantities $\widehat W_{i,\le c}(t)$ are (expected) total
workloads, since they incorporate work currently in the system
which will be routed to buffer $i$, regardless where the work resides in the network at time $t$.

\paragraph{Three types of jobs \& weights.}
We distinguish three types of jobs:
\begin{itemize}
\item[]
\begin{itemize}
\item[Type 1.] Jobs with counter $>D$
\item[Type 2.] Jobs on a path of
length $=D$ in $\cP$.
\item[Type 3.] Jobs on a path of
length $<D$ in $\cP$.
\end{itemize}
\end{itemize}
Hence, the counters of jobs of Type 2 and 3 cannot exceed $D$. We
further note that each job (regardless whether it
is currently being processed or not) can
compute its expected total remaining service requirement
in the future (under the network process $\{Y(t)\}$), and we call this quantity
the weight of the job.
Namely, if a job is not currently being processed, its weight is
\begin{eqnarray*}
&&\sum_{i\in \I} m_i\,E\big[\mbox{the number of times it visits
buffer $i$ before it leaves the network}\big],
\end{eqnarray*} which includes the current buffer of the job.
Thus, the weight of a waiting job, not currently being processed,
can be calculated as follows.
\begin{itemize}
\item[]
\begin{itemize}
\item[Type 1.] The weight of a waiting job of Type 1 in buffer
$i\in\I$ is
$${\bf e}_i^T(I+P+P^2+\cdots)\bm~=~{\bf e}_i^T(I-P)^{-1}\bm,$$
where $\bm=[m_i]\in\mathbb{R}_+^I$ and ${\bf e}_i\in \{0,1\}^I$ is
the unit vector with zeros except for its $i$-th coordinate which
equals to one (both are column vectors).
\item[Type 2.] The weight of a waiting job of Type 2 in
buffer $i_{nm}$ on a path $P_n$ of length $D$ (i.e., $|P_n|=D$) is
$${\bf e}_{i_{nD}}^T(P+P^2+P^3+\cdots)\bm+\sum_{k=m}^{D} m_{i_{nk}}
~=~{\bf e}_{i_{nD}}^TP(I-P)^{-1}\bm+\sum_{k=m}^{D} m_{i_{nk}},$$
where we note that $i_{nD}$ is the last buffer on $P_n$.
\item[Type 3.] The weight of a waiting job of Type 3 in
buffer $i_{nm}$ on a path $P_n$ of length $<D$ is
$$\sum_{k=m}^{|P_n|} m_{i_{nk}}.$$
\end{itemize}
\end{itemize}
On the other hand, for jobs currently being processed, each weight
is calculated as
\begin{eqnarray*}
&&\mbox{the current remaining service requirement} \\
&&\qquad+\sum_{i\in \I} m_i\,E\big[\mbox{the number of times it
visits buffer $i$ before it leaves the network}\big],
\end{eqnarray*}
where the latter number does not include the current buffer of the job. Now
let $M_1(t)$, $M_2(t)$ and $M_3(t)$ be the total weights of jobs of
Type 1, 2 and 3 at time $t$, respectively. These total
weights are fully determined by the network state information
$Y(t)$:
\begin{eqnarray*}
M_1(t)&=&\sum_{i\in\I,c>D} Q_{i,c}(t)\,{\bf e}_i^T(I-P)^{-1}\bm
+\sum_{i\in\I,c>D,j\in\J}
{\bf 1}^+_{
V^j_{i,c}(t)}\left[V^j_{i,c}(t)+
{\bf e}_i^TP(I-P)^{-1}\bm\right],\\
M_2(t)&=&\sum_{n\leq (I+1)^D,m\leq D: |P_n|=D}Q_{nm}(t)\left[{\bf
e}_{i_{nD}}^TP(I-P)^{-1}\bm+\sum_{k=m}^{D} m_{i_{nk}}\right]\\
&&\quad+\sum_{n\leq (I+1)^D,m\leq D,j\in J: |P_n|=D}
{\bf 1}^+_{
V^j_{nm}(t)}\left[V^j_{nm}(t)+{\bf e}_{i_{nD}}^TP(I-P)^{-1}\bm+
\sum_{k=m+1}^{D}
m_{i_{nk}}\right],\\
M_3(t)&=&\sum_{n\leq (I+1)^D,m\leq D:
|P_n|<D}\left[Q_{nm}(t)\sum_{k=m}^{|P_n|} m_{i_{nk}}\right]\\
&&\qquad+\sum_{n\leq (I+1)^D,m\leq D,j\in\J: |P_n|<D}
{\bf 1}^+_{
V^j_{nm}(t)}\left[V^j_{nm}(t)+\sum_{k=m+1}^{|P_n|}
m_{i_{nk}}\right],\end{eqnarray*}
where we recall that ${\bf 1}^+_x=1$ if $x>0$, and ${\bf 1}^+_x=0$ otherwise.

\subsubsection{Three
Auxiliary Lemmas and Proof of Proposition \ref{lem2}}\label{sec:skeleton} We
state and prove the following three key lemmas, which we prove in Section \ref{sec:pflemkey0}, \ref{sec:pflemkey1} and
\ref{sec:pflemkey2}, respectively.
To simplify notation, we will use $\mathcal C$ to denote a finite
constant which only depends on the matrix $Z$ given in Theorem
\ref{thm:main} or on the predefined network parameters from Section
\ref{sec:defSPN}. Its precise value can be different from line to
line.
\begin{lemma}\label{lem:key0}
There exists a constant $\mathcal C<\infty$ such that
 for all $t\geq 0$, 
\begin{eqnarray*}
E\left[\|\bV(t+1)\|_2^2~|~\cF(t)\right]&\leq&\|\bV(t)\|_2^2
 -\beta_{\min}\|\bV(t)\|_1 + \mathcal C,
 \end{eqnarray*}
 where $\beta_{\min}=\min_{j\in\J}\beta_j$.
\end{lemma}
\begin{lemma}\label{lem:key1}
	Suppose that there exists a symmetric matrix
	$Z\in\mathbb{R}_+^{I\times I}$ such that $\Lloc(x) = x^T Z x$
	is a local Lyapunov function with slack
	$\varepsilon\geq0$, and that either Condition C1 or {C2$\,^\prime$ from}
	Theorem \ref{thm:main} holds.
	Then given $D\in\mathbb N$, there exist
constants $\upsilon, \mathcal C=\cC(D)\in(0,\infty)$ such that for
all $c\leq D$, $t\geq 0$, 
\begin{eqnarray*}
\lefteqn{
E\left[\mathcal L_{\text{\em local}}\left(\widehat{\bW}_{\leq
c}(t+1)\right)~\Big|~\cF(t)\right]}\\&\leq&\mathcal L_{\text{\em
local}}\left(\widehat{\bW}_{\leq c}(t)\right)- \upsilon
\|{\bQ}_{\leq c}(t)\|_1+\cC\left(\|\bV(t)\|_1+\|{\bQ}_{<
c}(t)\|_1+1\right).
\end{eqnarray*}
\end{lemma}
\begin{lemma}\label{lem:key2}
	
Consider $\varepsilon>0$.
If the network is synchronized or if Condition C2 from Theorem \ref{thm:main} holds,
then there exist
constants $D, T\in\mathbb N$ and $\cC,\gamma_1,\gamma_2\in(0,\infty)$ such that  for
all $t\geq 0$, 
$$E\left[\mathcal G(Y(t+T))^2~|~\cF(t)\right]~\leq~
\mathcal G(Y(t))^2 -\gamma_2 \|\bQ_{>D}(t)\|_1+\cC\left(\|\bQ_{\leq
D}(t)\|_1+\|\bV(t)\|_1+1\right),$$
where $\mathcal G(Y(t)):=M_1(t)+M_2(t)+\gamma_1\|\bV(t)\|_1$.
\end{lemma}
Lemma \ref{lem:key2} is not needed for the proof of Proposition
\ref{lem2} if all routes are bounded.
Hence, for networks with bounded routes, $\varepsilon=0$ is allowed and
only Condition {\em C2$\,^\prime$} is needed for the desired stability.
The right-hand sides of the inequalities in Lemmas \ref{lem:key0} --
\ref{lem:key2} provide `negative drifts' on $\|\bV(t)\|_1$,
$\|{\bQ}_{\leq D}(t)\|_1$ and $\|{\bQ}_{>D}(t)\|_1$, respectively.
By appropriately weighing the functions in these three lemmas, we
shall construct an appropriate
(global) Lyapunov function $\Lglo$. To this end, %
from Lemma \ref{lem:key0}, we obtain
\begin{eqnarray}
	E\left[\|\bV(t+T)\|_2^2~|~\cF(t)\right]&\leq&
	E\left[\|\bV(t+T-1)\|_2^2~|~\cF(t)\right]
	 -\beta_{\min}E\left[\|\bV(t+T-1)\|_1~|~\cF(t)\right] + \mathcal C,\notag\\	
	&\leq&
	E\left[\|\bV(t+T-1)\|_2^2~|~\cF(t)\right]+\cC\notag\\
	&\leq&
	E\left[\|\bV(t+T-2)\|_2^2~|~\cF(t)\right]+2\cC\notag\\
	&&\qquad\qquad\dots\notag\\
&\leq&E\left[\|\bV(t+1)\|_2^2~|~\cF(t)\right]+(T-1)\cC\notag\\
&\leq&\|\bV(t)\|_2^2
 -\beta_{\min}\|\bV(t)\|_1 + T\mathcal C,\label{eq:bVbound1}
 \end{eqnarray}
where 
$T$ is the constant from Lemma \ref{lem:key2}.
We write this inequality as
\begin{eqnarray}
&&E\left[\frac1{\beta_{\min}}\|\bV(t+T)\|_2^2\,\Big|\,\cF(t)\right]~\leq~\frac1{\beta_{\min}}\|\bV(t)\|_2^2
 -\|\bV(t)\|_1 + \mathcal C,\label{eq1:pfmainthm}
\end{eqnarray}
where the constant $\cC$ is redefined appropriately.
We now argue that, similarly, Lemma~\ref{lem:key1} implies that
\begin{eqnarray}
&& E\left[\Z\left(\widehat{\bW}_{\leq
c}(t+T)\right)\,\Big|\,\cF(t)\right]\notag\\
&&\qquad\qquad\qquad\leq~\Z\left(\widehat{\bW}_{\leq
c}(t)\right)- \upsilon \|{\bQ}_{\leq
c}(t)\|_1+\cC\left(\|\bV(t)\|_1+\|{\bQ}_{<
c}(t)\|_1+1\right),\label{eq2:pfmainthm}
\end{eqnarray}
where $\cC$ is some (large) constant which may different from the one in Lemma \ref{lem:key1}.
To see why this holds, we use the same argument that led to (\ref{eq:bVbound1})
with the additional observation that, for some constant $\cC'$,
\begin{eqnarray}
	\sup_{s\in\{0,1,\ldots,T\}} \left(E\left[\|\bV(t+s)\|_1~|~\cF(t)\right]-\|\bV(t)\|_1\right)&\le&\cC',\label{eq:bVbound2}\\
	\sup_{s\in\{0,1,\ldots,T\}} \left(E\left[\|\bQ(t+s)\|_1~|~\cF(t)\right]-\|\bQ(t)\|_1\right)&\le&\cC'.\label{eq:bQbound}
\end{eqnarray}
Here \eqref{eq:bVbound2} can be derived along the lines of \eqref{eq:bVbound1},
and \eqref{eq:bQbound} follows from \eqref{eq:arrivalassum}
and the observation that the change in queue length is majorized by the number of external job arrivals.

We now show that Proposition \ref{lem2} follows from these three
lemmas, where the last lemma is not needed if all routes are bounded.
We consider the following Lyapunov function $\Lglo$: 
\begin{itemize}
	\item If all routes are bounded,
	\begin{equation*}
		\Lglo(Y(t))~:=~\sum_{c=1}^D \left(\frac{\upsilon}{2\mathcal C}\right)^c \Z\left(\widehat{\bW}_{\leq c}(t)\right)
	+\frac{2\mathcal C}{\beta_{\min}}
	\,\|\bV(t)\|_2^2,\end{equation*}
	where we choose $D<\infty$ such that $P^D=0$, and
	$\upsilon$ comes from Lemma \ref{lem:key1}.
	\item Otherwise,
\begin{equation*}
	\Lglo(Y(t))~:=~\sum_{c=1}^D \left(\frac{\upsilon}{2\mathcal C}\right)^c \Z\left(\widehat{\bW}_{\leq c}(t)\right)
+\frac{2\mathcal C}{\beta_{\min}}
\,\|\bV(t)\|_2^2+\frac{\xi}{2\cC}\,\mathcal G(Y(t))^2,\end{equation*} where
$\upsilon$ and $D$ is from Lemma \ref{lem:key1} and \ref{lem:key2}, respectively,
$\xi:=\xi(\upsilon,\cC,D)={\upsilon}\left(\frac{\upsilon}{2\mathcal C}\right)^D$
and $\cC$ is a large enough constant chosen so that
it can be used for
Lemma \ref{lem:key2} as well as for \eqref{eq1:pfmainthm} and \eqref{eq2:pfmainthm}.
\end{itemize}
We focus on proving Proposition \ref{lem2} for the case of unbounded routes,
but all arguments go through for the other case.
Without loss of generality, we assume that $$\xi,\upsilon,\gamma_2<1<\cC.$$
The property
$\sup_{Y\in\Omega_Y^*}{\Lglo(Y)}/{\mathcal |Y|^2}<\infty$ in
Proposition \ref{lem2} is readily seen to hold. 
To derive the negative drift property, we observe
that Lemma \ref{lem:key2} in conjunction with \eqref{eq1:pfmainthm} and \eqref{eq2:pfmainthm} imply that
\begin{eqnarray*}
&&E\left[\Lglo(Y(t+T))-\Lglo(Y(t))\,\big|\,\cF(t)\right]\\
 &&\qquad\leq\cC^*-2\mathcal C
\|{\bV}(t)\|_1+\sum_{c=1}^D \left(\frac{\upsilon}{2\mathcal C}\right)^c
\Big[- \upsilon \|{\bQ}_{\leq c}(t)\|_1+\mathcal C \|{\bQ}_{<
c}(t)\|_1+\mathcal
C \|{\bV}(t)\|_1\Big]\\
&&\qquad\qquad
-\frac{\gamma_2\xi}{2\cC} \|\bQ_{>D}(t)\|_1+\frac{\xi}{2}\left(\|\bQ_{\leq D}(t)\|_1+\|\bV(t)\|_1+1\right)\\
&&\qquad\leq\cC^*-\mathcal C \|{\bV}(t)\|_1+\sum_{c=1}^D
\left(\frac{\upsilon}{2\mathcal C}\right)^c \Big[- \upsilon
\|{\bQ}_{\leq c}(t)\|_1+\mathcal C \|{\bQ}_{< c}(t)\|_1\Big]\\
&&\qquad\qquad
-\frac{\gamma_2\xi}{2\cC} \|\bQ_{>D}(t)\|_1+\frac{\xi}{2}\left(\|\bQ_{\leq D}(t)\|_1+\|\bV(t)\|_1+1\right),
\end{eqnarray*} where $\cC^*$ is some (large enough) constant and we use $\upsilon\leq \cC$
for the last inequality.
The sum in this expression can be bounded as follows:
\begin{eqnarray*}
\lefteqn{\sum_{c=1}^D \left(\frac{\upsilon}{2\mathcal C}\right)^c \Big[- \upsilon
	\|{\bQ}_{\leq c}(t)\|_1+\mathcal C \|{\bQ}_{< c}(t)\|_1\Big]}\\	
&=& -\upsilon\sum_{c=1}^D
\left(\frac{\upsilon}{2\mathcal C}\right)^c \|{\bQ}_{\leq c}(t)\|_1
+\frac{\upsilon}2\sum_{c=1}^D
\left(\frac{\upsilon}{2\mathcal C}\right)^{c-1} \|{\bQ}_{< c}(t)\|_1\\
&=&-{\upsilon}\left(\frac{\upsilon}{2\cC}\right)^D\|{\bQ}_{\leq D}(t)\|_1-\left(
\upsilon-\frac{\upsilon}2\right)\sum_{c=1}^{D-1}
\left(\frac{\upsilon}{2\mathcal C}\right)^c \|{\bQ}_{\leq c}(t)\|_1\\
&\leq&- {\upsilon}\left(\frac{\upsilon}{2\cC}\right)^D\|{\bQ}_{\leq D}(t)\|_1
~=~-\xi\|{\bQ}_{\leq D}(t)\|_1.
\end{eqnarray*}
After combining the preceding two displays, we obtain
the desired negative drift property:
\begin{eqnarray*}
	\lefteqn{E\left[\Lglo(Y(t+T))-\Lglo(Y(t))\,\big|\,\cF(t)\right]}\\
&\le&\cC^*-\mathcal C \|{\bV}(t)\|_1-\xi\|{\bQ}_{\leq D}(t)\|_1
-\frac{\gamma_2\xi}{2\cC} \|\bQ_{>D}(t)\|_1+\frac{\xi}{2}\left(\|\bQ_{\leq D}(t)\|_1+\|\bV(t)\|_1+1\right)\\
&\leq&\cC^*+\frac{\xi}2-\frac{\gamma_2\xi}{2\cC}\left(\|{\bV}(t)\|_1+
\|{\bQ_{\leq D}}(t)\|_1+\|{\bQ_{>D}}(t)\|_1\right)\\
&=&\cC^*+\frac{\xi}2-\frac{\gamma_2\xi}{2\cC}\left(\|{\bV}(t)\|_1+\|{\bQ}(t)\|_1\right),
\end{eqnarray*} where we use $\xi, \upsilon,\gamma_2<1$ and $\cC>1$. This completes the proof
of Proposition \ref{lem2}.

\subsubsection{Proof of Lemma \ref{lem:key0}}\label{sec:pflemkey0}
Recall that $\|\bV(t)\|_2^2$ is the sum of squares of the remaining
service requirements of jobs being processed by some activity at
time $t$, i.e., $$\|\bV(t)\|_2^2~=~\sum_{j\in\J}
V^j(t)^2\qquad\mbox{and}\qquad\|\bV(t)\|_1~=~\sum_{j\in\J} V^j(t).$$
On the event $\{V^{j}(t)\leq \beta_j\}$, activity $j$ has to restart
before time $t+1$ and hence,
\begin{eqnarray*}
	\lefteqn{E\left[V^{j}(t+1)^2~|~\cF(t)\right]}\\&\leq&
	E\left[V^{j}(t+1)^2~|~V^{j}(t+1)>0,\cF(t)\right]\\
	&\leq&E\left[\Gamma_{i_j,t+1}^2~\big|~V^{j}(t+1)>0,\cF(t)\right]\\
	&=& E\left[E\left[\Gamma_{i_j,t+1}^2~\big|~V^{j}(t+1)>0,\cF(t),J^{(t+1)}
\right]~\big|~V^{j}(t+1)>0,\cF(t)\right],
\end{eqnarray*}
where we let $\Gamma_{i_j,t+1}$ 
denote the service time generated by the job
being processed by activity $j$ at time $t+1$
and we write $$J^{(t+1)}:=t+1-\frac{\Gamma_{i_j,t+1}-V^j(t+1)}{\beta_j}\in[t,t+1]$$
for the time when this job starts its service. Note that, again on
the event $\{V^j(t)\le \beta_j\}$,
\begin{eqnarray}
\lefteqn{E\left[\Gamma_{i_j,t+1}^2~\big|~V^{j}(t+1)>0,\cF(t),J^{(t+1)}\right]}\notag\\
&=&\int^\infty_0 \Pr\left[\Gamma_{i_j,t+1}^2>x~\big|~V^{j}(t+1)>0,\cF(t),J^{(t+1)}\right]~dx\notag\\
&=&\int^\infty_0 \Pr\left[\Gamma_{i_j}^2>x~\big|\Gamma_{i_j}>\beta_j\left(t+1-J^{(t+1)}\right)\right]~dx\notag\\
&=&\int^\infty_0 \frac{\Pr\left[\Gamma_{i_j}^2>x,\Gamma_{i_j}>\beta_j\left(t+1-J^{(t+1)}\right)\right]}{\Pr\left[\Gamma_{i_j}
>\beta_j\left(t+1-J^{(t+1)}\right)\right]}~dx\notag\\
&\leq&\int^\infty_0 \frac{\Pr\left[\Gamma_{i_j}^2>x\right]}{\Pr[\Gamma_{i_j}>\beta_j]}~dx\notag\\
&=&\frac{E[\Gamma_{i_j}^2]}{\Pr[\Gamma_{i_j}>\beta_j]},\notag
	\end{eqnarray}
where we recall that $\Gamma_{i_j}$ stands for a generic service time for buffer $i_j$.
We have thus established that on $\{V^{j}(t)\leq \beta_j\}$,
\begin{eqnarray}\label{eq1:lemkey0}
	E\left[V^{j}(t+1)^2~|~\cF(t)\right]&\leq&
	\begin{cases}
		\frac{E[\Gamma_{i_j}^2]}{\Pr[\Gamma_{i_j}>\beta_j]}&\mbox{if}~\Pr[\Gamma_{i_j}>\beta_j]>0\\
		~~~~\quad(\beta_j)^2&\mbox{otherwise}
	\end{cases}.
\end{eqnarray}

On the other hand, on the event $\{V^{j}(t)>\beta_j\}$,
\begin{equation}\label{eq2:lemkey0}
	V^{j}(t+1)^2~=~(V^{j}(t)-\beta_j)^2~\le~V^{j}(t)^2-\beta_jV^{j}(t)+\beta_j^2.
\end{equation}
Hence, combining \eqref{eq1:lemkey0} and \eqref{eq2:lemkey0}, we find that for some constant $\cC<\infty$,
$$E[V^{j}(t+1)^2~|~\cF(t)]~\leq~
V^{j}(t)^2-\beta_jV^{j}(t)+\cC,$$ which leads to the desired
conclusion
of Lemma \ref{lem:key0}. 

\subsubsection{Proof of Lemma \ref{lem:key1}}\label{sec:pflemkey1}
For notational convenience, we stick to the case $t=0$ in the
conclusion of Lemma \ref{lem:key1}. Namely, we show that
\begin{eqnarray*}
E\left[\Z\left(\widehat{\bW}_{\leq
c}(1)\right)~\big|~\cF(0)\right]&\leq&\Z\left(\widehat{\bW}_{\leq
c}(0)\right)- \upsilon \|{\bQ}_{\leq
c}(0)\|_1+\cC\left(\|\bV(0)\|_1+\|{\bQ}_{<
c}(0)\|_1+1\right).
\end{eqnarray*}
However, all arguments go through for general $t>0$.

\paragraph{Non-synchronized network.}
We first consider the case when the network may not be synchronized.
The first step in the proof is the observation that the
schedule $\bsigma(t)$ under the $\varepsilon$-LRFS policy is always
maximal with respect to the vector $\bw=[w_i]=\left[m_i Q_{i,\leq
c}(t)\right]$ for any $c\geq 1$. Consequently, since
the local quadratic Lyapunov function
$\Lloc(x)=x^TZx$ satisfies \eqref{eq:localcondition}, we obtain that
for $t\in [0,1]$,
\begin{equation*}
2\sum_{(i,\ell)\in \I\times\I} Z_{i\ell}\,w_\ell
\left(\rho_i+\varepsilon m_i-\sum_{j\in \J_i}\beta_j\sigma_j(t)
\right)~ \leq~ -\eta \|\bw\|_1+\cC,\end{equation*} where we further
use the observation that
\begin{eqnarray}
	&& \left(w_i+\rho_i+
	\varepsilon m_i-\sum_{j\in
\J_i}\beta_j\sigma_j(t)\right)
\left(w_\ell+\rho_\ell+\varepsilon m_\ell-\sum_{j\in
\J_\ell}\beta_j\sigma_j(t)\right)-w_i w_\ell\notag\\
&&\qquad\leq~w_\ell\left(\rho_i+\varepsilon m_i-\sum_{j\in
\J_i}\beta_j\sigma_j(t)\right)+w_i\left(\rho_\ell+\varepsilon m_\ell-\sum_{j\in
\J_\ell}\beta_j\sigma_j(t)\right)+\cC.\label{eq0-2:lemkey1}
\end{eqnarray}
We remind the reader that we write $\cC$ for a finite
constant which may differ from line to line. Since we set $w_i=m_i
Q_{i,\leq c}(t)$, it follows that
\begin{equation*}
2\sum_{(i,\ell)\in \I\times\I} Z_{i\ell}\,m_\ell Q_{\ell,\leq c}(t)
\left(\rho_i+\varepsilon m_i-\sum_{j\in \J_i}\beta_j\sigma_j(t)
\right)~ \leq~ -\eta \sum_{i\in I} m_i Q_{i,\leq
c}(t)+\cC.\end{equation*}
After taking conditional expectations given $\cF(0)$ on both sides in the above inequality,
we obtain
\begin{equation*}
2\sum_{(i,\ell)\in \I\times\I} Z_{i\ell}\,m_\ell\cdot E\left[{Q}_{\ell,\leq
c}(t)  \left(\rho_i+\varepsilon m_i-\sum_{j\in
\J_i}\beta_j\sigma_j(t)\right)~\Bigg|~\cF(0)\right]~ \leq~ -\eta \sum_{i\in I}
m_i Q_{i,\leq c}(0)+\cC,
\end{equation*}
where we use that
$E[{Q}_{\ell,\leq
c}(t)~|~\cF(0)]\geq{Q}_{\ell,\leq c}(0)-\cC'$ for $t\in [0,1]$ for some constant $\cC'$.
This can be verified by suppressing any arrivals and letting all activities work,
and then using the standard fact from renewal theory that any renewal function is finite.
Similarly, we obtain from $E\left[{Q}_{\ell,\leq c}(t)~\big|~\sum_{j\in
\J_i}\beta_j\sigma_j(t),\cF(0)\right]\geq {Q}_{\ell,\leq c}(0)-\cC'$ for $t\in[0,1]$ that
\begin{equation}\label{eq0:lemkey1}
2\sum_{(i,\ell)\in \I\times\I} Z_{i\ell}\,m_\ell {Q}_{\ell,\leq
c}(0) \cdot E\left[\rho_i+\varepsilon m_i-\sum_{j\in
\J_i}\beta_j\sigma_j(t)~\Bigg|~\cF(0)\right]~ \leq~ -\eta \sum_{i\in I}
m_i Q_{i,\leq c}(0)+\cC,
\end{equation}
where the constant $\cC$ again has to be redefined appropriately.
We leave this inequality for later use.

The second step in the proof is to bound
$E\left[\widehat{W}_{i,\leq c}(1)-\widehat{W}_{i,\leq
c}(0)~\Big|~\cF(0)\right]$ for fixed $i$.
Let $\widehat{A}_{i,\leq c}$ be the
number of job arrivals contributing to an increase in $\widehat{Q}_{i,\leq
c}(\cdot)$ during the time interval $[0,1)$. Then, one can check that
for $\widehat{\bA}_{\leq c}=\left[\widehat{A}_{i,\leq c}\right]$,
\begin{equation*}
E\left[\widehat{\bA}_{\leq c}~\big|~\cF(0)\right]~\leq~
\left(I+P+\dots
+P^c\right)\balpha~\leq~\blambda,
\end{equation*}
since we assume \eqref{eq:arrivalassum}.
Define the following quantities.
\begin{itemize}
\item ${D}_{i,\leq c}$ and ${D}_{i,>c}$ are the numbers of jobs in buffer $i$
which start their service during the time interval $[0,1)$ and with
counter $\leq c$ and $>c$, respectively.
\item $R_{i,\leq c}$ and $R_{i,> c}$ are the total amounts of service
times generated by jobs contributing to ${D}_{i,\leq c}$ and
${D}_{i,>c}$, respectively, again during the time interval $[0,1)$.
We stress that the contribution of each job to these quantities may exceed the service time
it receives during the interval $[0,1)$.
\item $R_{i,> c}^{(L)}$ are the total amounts of
service times generated by jobs contributing to $D_{i,>c}$ due to
the LRFS policy (i.e., step 4 
in Definition \ref{def:epLRFS}). In addition, we set $R_{i,> c}^{(M)}=R_{i,> c}-R_{i,>c}^{(L)}$.
\end{itemize} Then, we have that
\begin{equation}
\Delta_{i,\leq c}~:=~\widehat{W}_{i,\leq c}(1)-\widehat{W}_{i,\leq
c}(0)
~=~m_i\left(\widehat{A}_{i,\leq c} - {D}_{i,\leq c}\right)+R_{i,\leq
c} +R_{i,> c}- \sum_{j\in \J_i}\int^1_0 \beta_j\sigma_j(t)\, dt.\label{eq:defDelta}
\end{equation}
Taking conditional expectations given $\cF(0)$ on both sides, we
obtain\begin{eqnarray}&& E\left[\Delta_{i,\leq c}~|~\cF(0)\right]\notag\\
&&\qquad\leq
m_i\lambda_i +E\left[- m_i{D}_{i,\leq c}+R_{i,\leq
c}~|~\cF(0)\right]+E\left[R_{i,> c}~|~\cF(0)\right] - E\left[\sum_{j\in
\J_i}\int^1_0 \beta_j\sigma_j(t)\, dt~\Bigg|~\cF(0)\right]\notag\\
&&\qquad=\rho_i+E\left[R_{i,> c}~|~\cF(0)\right]  -E\left[\sum_{j\in
\J_i}\int^1_0 \beta_j\sigma_j(t)\,
dt~\Bigg|~\cF(0)\right].\label{eq3:lemkey1}
\end{eqnarray}

Now we bound $E\left[R_{i,> c}~|~\cF(0)\right]$ in the above
inequality, or equivalently $E\left[R_{i,>
c}^{(L)}~\Big|~\cF(0)\right]+E\left[R_{i,> c}^{(M)}~\Big|~\cF(0)\right]$
since $R_{i,>c}=R_{i,> c}^{(M)}+R_{i,> c}^{(L)}$. First, one can
check that
\begin{equation*}
	E\left[R_{i,> c}^{(M)}~\Big|~\cF(0)\right]~\leq~\varepsilon m_i.
	\end{equation*}
This is because the expected number of jobs in buffer $i$
which start their service during the time interval $[0,1)$ due to
step 3-1 of the $\varepsilon$-LRFS policy in Definition \ref{def:epLRFS}
is at most $\varepsilon$ since each timer $\mathcal T^{(h)}(t)$
is zero at most once during this time
interval.
On the
other hand, to bound $E\left[R_{i,> c}^{(L)}~\Big|~\cF(0)\right]$,
consider activity-interchangeable buffers $\ell$ and $i$ (which
includes $\ell=i$). Let $\mathcal E_\ell$ be the event that every job in
the queue $Q_{\ell,\leq c}(0)$ (i.e., jobs with counter $\leq c$
waiting in buffer $\ell$ at time $0$)
starts service 
during the time interval $[0,1]$. One can observe that $R_{i,> c}^{(L)}=0$
on the complementary event $\overline{\mathcal E}_\ell$. 
We let
$$N:~=~\sum_{w=1}^{Q_{\ell,\leq c}(0)-|\J_\ell|} \Gamma_{\ell,w},$$ where
$\{\Gamma_{\ell,w}\}$ are identical random variables with mean $m_\ell$
and variance $\varsigma_\ell^2<\infty$. We first consider the case when
\begin{equation}\label{eq:assumptionN}
	E[N~|~Q_{\ell,\leq c}(0)]=m_\ell
(Q_{\ell,\leq c}(0)-|\J_\ell|)>\beta_{\max}|\J_\ell|.\end{equation} Since $\mathcal E_\ell$
occurs only if at least $Q_{\ell,\leq c}(0)-|\J_\ell|$ jobs complete
their service requirements, we obtain the following on the event that \eqref{eq:assumptionN} holds:
\begin{eqnarray}
\Pr[\mathcal E_\ell~|~\cF(0)]&\leq&\Pr\left[N\leq \beta_{\max}|\J_\ell|~|~\cF(0)\right]\notag\\
 &=&\Pr\left[N-E[N~|~Q_{\ell,\leq c}(0)]\leq
\beta_{\max}|\J_\ell|-E[N~|~Q_{\ell,\leq c}(0)]~|~\cF(0)\right]\notag\\
 &\leq&\Pr\left[\big(N-E[N~|~Q_{\ell,\leq c}(0)]\big)^2\geq
\big(E[N~|~Q_{\ell,\leq c}(0)]-\beta_{\max}|\J_\ell|\big)^2~|~\cF(0)\right]\notag\\
&\leq&\frac{\varsigma_\ell^2 (Q_{\ell,\leq
c}(0)-|\J_\ell|)}{\big(E[N~|~Q_{\ell,\leq c}(0)]-\beta_{\max}|\J_\ell|\big)^2}\notag\\
&=&\frac{\varsigma_\ell^2 (Q_{\ell,\leq
c}(0)-|\J_\ell|)}{\big(m_\ell(Q_{\ell,\leq c}(0)-|\J_\ell|)-\beta_{\max}|\J_\ell|\big)^2}\notag\\
&\leq&\frac{\cC}{Q_{\ell,\leq c}(0)+1},\label{eq1:lemkey1}
\end{eqnarray}
where $\beta_{\max}=\max_{j\in \J}\beta_j$, $\cC$ is some (finite)
constant depending on $\varsigma_\ell, m_\ell,|\J_\ell|$ and we use
Markov's inequality in conjunction with \eqref{eq:assumptionN}. On the event that \eqref{eq:assumptionN}
does not hold, i.e., when $Q_{\ell,\leq c}(0)$ is
bounded above by $|\J_\ell|\left(\beta_{\max}+\frac1{m_\ell}\right)$,
one can redefine the constant $\cC$ so that
\eqref{eq1:lemkey1} holds. Hence, \eqref{eq1:lemkey1} always holds.

Using \eqref{eq1:lemkey1}, it follows that
\begin{eqnarray*}
E[R_{i,> c}~|~\cF(0)]&=&E\left[R_{i,> c}^{(M)}~\Big|~\cF(0)\right]+E\left[R_{i,>
c}^{(L)}~\Big|~\cF(0)\right]\\
&\leq&\varepsilon m_i+ \Pr[\mathcal
E_\ell~|~\cF(0)]\cdot E\left[R_{i,> c}^{(L)}~\big|~\mathcal E_\ell,\cF(0)\right]\\
&\leq&\varepsilon m_i+ \frac{\cC}{Q_{\ell,\leq c}(0)+1}\cdot
E\left[R_{i,> c}^{(L)}~\big|~\mathcal E_\ell,\cF(0)\right]\\
&\leq&\varepsilon m_i+ \frac{\cC}{Q_{\ell,\leq c}(0)+1},
\end{eqnarray*}
where the last inequality
requires that the constant $\cC$ has to be redefined appropriately, since
$E\left[R_{i,> c}^{(L)}~\big|~\mathcal E_\ell, \cF(0)\right]\le \cC'$
as can be seen using arguments similar to those leading up to \eqref{eq1:lemkey0}.
Together with \eqref{eq3:lemkey1}, this leads to
\begin{eqnarray}
E\left[\Delta_{i,\leq c}~|~\cF(0)\right]&\leq& \rho_i+\varepsilon
m_i+\frac{\cC}{Q_{\ell,\leq c}(0)+1} -E\left[\sum_{j\in
\J_i}\int^1_0 \beta_j\sigma_j(t)\, dt~\Bigg|~\cF(0)\right],\label{eq4:lemkey1}
\end{eqnarray}
for any activity-interchangeable buffers $\ell$ and $i$.

The third step in the proof for the non-synchronized case is to prove the conclusion of Lemma
\ref{lem:key1}.
By a similar argument as in \eqref{eq0-2:lemkey1},
the claim follows with $\upsilon=\eta\min_{i\in\I}m_i$
after we show that
\begin{equation}\label{eq2:lemkey1}
2\sum_{(i,\ell)\in \I\times\I} Z_{i\ell}\,\widehat{W}_{\ell,\leq
c}(0)\,E[\Delta_{i,\leq c}~|~\cF(0)]~\leq~- \eta
\sum_{i\in\I}m_i{Q}_{i,\leq c}(0)+ \cC\sum_{i\in\I}
\Big[V_i(0)+{Q}_{i,< c}(0)+1\Big],
\end{equation}
where we use that
\begin{equation}\label{eq2-1:lemkey1}
	E\left[\Delta_{i,\leq c}\Delta_{\ell,\leq
c}~|~\cF(0)\right]~\leq~ \cC\end{equation}
for all $i,\ell$ and some constant $\cC$. To see that \eqref{eq2-1:lemkey1} holds,
it suffices to show that $E\left[\Delta_{i,\leq c}^2~|~\cF(0)\right]\leq \cC$
by the Cauchy-Schwarz inequality. This can be
shown using \eqref{eq:defDelta} and
$$E\left[\widehat{A}_{i,\leq c}^2~\big|~\cF(0)\right], E\left[D_{i,\leq c}^2~|~\cF(0)\right],
E\left[R_i^2~|~\cF(0)\right]~\leq~ \cC$$ where these bounds can be derived using
arguments similar to those leading up to \eqref{eq1:lemkey0}.
Since $Z_{i\ell},
E[\Delta_{i,\leq c}~|~\cF(0)]\leq\cC$ and
$$\widehat{W}_{\ell,\leq c}(0)~\leq~ {W}_{\ell,\leq c}(0)+\sum_{i\in\I}m_i{Q}_{i,< c}(0)
~=~m_\ell{Q}_{\ell,\leq c}(0)+V_\ell(0)+\sum_{i\in\I}m_i{Q}_{i,< c}(0),$$
the inequality in \eqref{eq2:lemkey1} reduces to 
\begin{equation}\label{eq5:lemkey1} 2\sum_{(i,\ell)\in
\I\times\I} Z_{i\ell}\,m_\ell{Q}_{\ell,\leq c}(0)\,E[\Delta_{i,\leq
c}~|~\cF(0)]~\leq~-\eta \sum_{i\in\I}m_i{Q}_{i,\leq c}(0)+\cC.\end{equation}
We prove this using \eqref{eq0:lemkey1}
and \eqref{eq4:lemkey1} in conjuction with Condition {\em C2$\,^\prime$} as follows:
\begin{eqnarray*}
&&2\sum_{(i,\ell)\in \I\times\I} Z_{i\ell}\,m_\ell{Q}_{\ell,\leq
c}(0)\,E[\Delta_{i,\leq c}~|~\cF(0)]\\
&&\qquad\leq~ 2\sum_{(i,\ell)\in \I\times\I}
Z_{i\ell}\,m_\ell{Q}_{\ell,\leq c}(0)\left(\rho_i+\varepsilon
m_i+\frac{\cC}{Q_{\ell,\leq c}(0)+1} -E\left[\sum_{j\in
\J_i}\int^1_0
\beta_j\sigma_j(t)\, dt~\Bigg|~\cF(0)\right]\right)\\
&&\qquad\leq~ 2\sum_{(i,\ell)\in \I\times\I}
Z_{i\ell}\,m_\ell {Q}_{\ell,\leq c}(0)\left(\rho_i+\varepsilon m_i
-E\left[\sum_{j\in \J_i}\int^1_0
\beta_j\sigma_j(t)\, dt~\Bigg|~\cF(0)\right]\right)+\cC\\
&&\qquad\leq~-\eta \sum_{i\in\I}m_i{Q}_{i,\leq
c}(0)+\cC,
\end{eqnarray*}
where we again remind the reader that the constant $\cC$ may differ
from line to line. This completes the proof of Lemma \ref{lem:key1}
for non-synchronized networks.
\paragraph{Synchronized network.}
Now we consider the case when the network is synchronized, i.e., Condition {\em C1}.
We establish the same three steps
as in the non-synchronized case.
In the non-synchronized case,
we used the fact that the schedule $\bsigma(t)$ under the
$\varepsilon$-LRFS policy is maximal with respect to $[m_iQ_{i,\leq c}(t)]$,
for which we required Condition {\em C2$\,^\prime$}.
In synchronized networks, as a first step in the proof,
we use a different (stronger) maximality property,
which allows us to relax Condition {\em C2$\,^\prime$} to Condition {\em C1}. To this end, we introduce some necessary notation.
We let $\sigma_{j,\leq c}(t)=1$ if activity $j$ processes a job with counter $\leq c$
at time $t$ (and $\sigma_{j,\leq c}(t)=0$ otherwise). Since $\sigma_{j,\leq c}(t)=
\sigma_{j,\leq c}(\lfloor t\rfloor)$ in synchronized networks, we write
$\sigma_{j,\leq c}=\sigma_{j,\leq c}(0)=\sigma_{j,\leq c}(t)$ for $t\in[0,1)$.
The main maximality property we use 
in synchronized networks is that, under the LRFS policy,
the schedule $\left[\sigma_{j,\leq c} I_{\I^{(h)}}(i_j)\right]$
is maximal with respect to $\left[Q_{i,\leq c}(0)I_{\I^{(h)}}(i)\right]$ for
each component $\I^{(h)}$.
Together with \eqref{eq:localcondition}, this implies that for every partition $\I^{(h)}$,
\begin{equation}\label{eq6-00:lemkey1}
2\sum_{(i,\ell)\in \I^{(h)}\times\I^{(h)}} Z_{i\ell}\,{Q}_{\ell,\leq
c}(0)\,\left(\rho_i+\varepsilon m_i-\sum_{j\in \J_i}
E\left[\sigma_{j,\leq c}~|~\cF(0),\mathcal E_{\mbox{\tiny LRFS}}^{(h)}\right]\right)~\leq~-\eta
\sum_{i\in\I}{Q}_{i,\leq c}(0)+\cC,
\end{equation}
where we let $\mathcal E_{\mbox{\tiny LRFS}}^{(h)}$
denote the event that at time $0$ the $\varepsilon$-LRFS policy for component $\I^{(h)}$
does not select a job for processing
through step 3-1 (see its description in Definition \ref{def:epLRFS}) and
all selected jobs are due the LRFS policy in
step 4. We stress that in synchronized networks, every processor completes
the service requirement of the job it processes at every integer time $t\in \mathbb Z_+$ and hence
$\Pr\left[\mathcal
E_{\mbox{\tiny LRFS}}^{(h_j)}\right]=1-\varepsilon$. Inequality \eqref{eq6-00:lemkey1}
is analogous to \eqref{eq0:lemkey1}, which concludes the first step in the non-synchronized case.

We proceed with the analog of the second step
from the non-synchronized case, i.e.,
bounding $E\left[\widehat{W}_{i,\leq c}(1)-\widehat{W}_{i,\leq
c}(0)~|~\cF(0)\right]$.
We stress that the definition of $\widehat{W}_{i,\leq c}(t)$
differs from the one used in non-synchronized networks, see Section \ref{sec:not}.
Scheduling decisions are only made at integer time epochs
(i.e., $\sigma(t)=\sigma(\lceil t\rceil)$) in synchronized
networks, so that
\begin{eqnarray}
	E\left[\widehat{W}_{i,\leq c}(1)-\widehat{W}_{i,\leq
	c}(0)~|~\cF(0)\right]&=&
	E\left[\widehat{A}_{i,\leq c}~|~\cF(0)\right]
	- \sum_{j\in \J_i} E\left[\sigma_{j,\leq c}~|~\cF(0)\right]\notag\\
&\leq&\rho_i	- \sum_{j\in \J_i} E\left[\sigma_{j,\leq c}~|~\cF(0)\right]\notag,
\end{eqnarray}
where we again let $\widehat{A}_{i,\leq c}$ be the
number of job arrivals contributing to an increase in $\widehat{Q}_{i,\leq
c}(\cdot)$ during the time interval $[0,1)$. On writing
$\Delta_{i,\leq c}:=\widehat{W}_{i,\leq c}(1)-\widehat{W}_{i,\leq
c}(0)$, we have
\begin{eqnarray}
	E\left[{\Delta}_{i,\leq c}~|~\cF(0)\right]
	&\leq&\rho_i	- \sum_{j\in \J_i} E\left[\sigma_{j,\leq c}~|~\cF(0)\right].
	\label{eq6-1:lemkey1}
\end{eqnarray}

The third step in the proof for the synchronized case is to prove the conclusion of Lemma
\ref{lem:key1}.
As in \eqref{eq5:lemkey1}, it suffices to prove that
$$2\sum_{(i,\ell)\in \I\times\I} m_\ell Z_{i\ell}\,{Q}_{\ell,\leq
c}(0)\,E[\Delta_{i,\leq c}~|~\cF(0)]~\leq~-\eta \sum_{i\in\I}m_i{Q}_{i,\leq
c}(0)+\cC. $$
Since $m_i=1$ and $V_i(0)\leq 1$ in synchronized networks, this reduces to
\begin{equation}\label{eq6:lemkey1} 2\sum_{(i,\ell)\in \I\times\I}
Z_{i\ell}\,{Q}_{\ell,\leq c}(0)\,E[\Delta_{i,\leq c}~|~\cF(0)]~\leq~-\eta
\sum_{i\in\I}{Q}_{i,\leq c}(0)+\cC.
\end{equation}
Combining \eqref{eq6-1:lemkey1} and \eqref{eq6:lemkey1}, it suffices to show that
\begin{equation}\label{eq7:lemkey1}
2\sum_{(i,\ell)\in \I\times\I} Z_{i\ell}\,{Q}_{\ell,\leq
c}(0)\,\left(\rho_i-\sum_{j\in \J_i}
E\left[\sigma_{j,\leq c}~|~\cF(0)\right]\right)~\leq~-\eta
\sum_{i\in\I}{Q}_{i,\leq c}(0)+\cC.
\end{equation}
For activity $j\in \J$,
writing $\I^{(h_j)}$ for the component of buffer $i_j$ (i.e., $i_j\in \I^{(h_j)}$), we have
\begin{eqnarray*}
	E\left[\sigma_{j,\leq c}~|~\cF(0)\right]&\geq&
	\Pr\left[\mathcal E_{\mbox{\tiny LRFS}}^{(h_j)}\right]
	E\left[\sigma_{j,\leq c}~|~\cF(0),\mathcal E_{\mbox{\tiny LRFS}}^{(h_j)}\right]\\
	&=&
(1-\varepsilon)	E\left[\sigma_{j,\leq c}~|~\cF(0),\mathcal E_{\mbox{\tiny LRFS}}^{(h_j)}\right]\\
	&\geq&	E\left[\sigma_{j,\leq c}~|~\cF(0),\mathcal E_{\mbox{\tiny LRFS}}^{(h_j)}\right]
-\varepsilon.
\end{eqnarray*}
Therefore, \eqref{eq7:lemkey1} follows after arguing that
\begin{equation*}
2\sum_{(i,\ell)\in \I\times\I} Z_{i\ell}\,{Q}_{\ell,\leq
c}(0)\,\left(\rho_i+\varepsilon m_i-\sum_{j\in \J_i}
E\left[\sigma_{j,\leq c}~|~\cF(0),\mathcal E_{\mbox{\tiny LRFS}}^{(h_j)}\right]\right)~\leq~-\eta
\sum_{i\in\I}{Q}_{i,\leq c}(0)+\cC.
\end{equation*}
The above inequality follows from \eqref{eq6-00:lemkey1} and Condition {\em C1}.
This completes the proof of Lemma \ref{lem:key1} for synchronized networks.

\subsubsection{Proof of Lemma \ref{lem:key2}}\label{sec:pflemkey2}
For notational convenience, we again restrict attention to the case $t=0$ in the
conclusion of Lemma \ref{lem:key2}, namely, we show that for some $D,T,\cC,\gamma_1, \gamma_2\in(0,\infty)$,
$$E\left[\mathcal G(Y(T))^2~|~\cF(0)\right]~\leq~
\mathcal G(Y(0))^2 -\gamma_2 \|\bQ_{>D}(0)\|_1+\cC\left(\|\bQ_{\leq
D}(0)\|_1+\|\bV(0)\|_1+1\right),$$
where we recall that $$\mathcal G(Y(t))=M_1(t)+M_2(t)+\gamma_1\|\bV(t)\|_1.$$
All arguments are
applicable for general $t\geq0$ as well. First observe that
$M(t):=M_1(t)+M_2(t)$ can only change through the following events
for jobs of Type 1 and 2.
\begin{itemize}
\item[] {\bf Arrivals.} $M(t)$ increases when new external arrivals of Type 2
occur. Note that there are no such external arrivals for
Type 1.
\item[] {\bf Routing.} $M(t)$ may increase or decrease when a job with counter $\geq D$ (i.e., Type 1 or 2) is
routed since the weight (i.e., future workload) of a job conditioned on
the buffer to which it has been routed
is different from the (unconditional) weight before it is routed.
However, $M(t)$ does not change when a job with counter $<D$
(i.e., Type 2) is routed since it is routed deterministically.
\item[] {\bf Starting Service.}
$M(t)$ may increase or decrease when a job
of Type 1 or 2 begins service, generating its service time at this point.
Assuming the job is served from buffer $i$, then
$M(t)$ increases if the random service time is larger than its mean $m_i$, and decreases
otherwise.
\item[] {\bf Being in Service.}
$M(t)$ decreases when a job of Type 1 or 2 is currently being
processed.
\end{itemize}
Now we express $M(t)$ as follows: for $t\in\mathbb{N}$,
$$M(t)~=~M(0)+M_{\text{arrival}}(t)+M_{\text{routing}}(t)+M_{\text{s-service}}(t)+M_{\text{b-service}}(t),$$
where $M_{\text{arrival}}(t)$, $M_{\text{routing}}(t)$,
$M_{\text{s-service}}(t)$ and $M_{\text{b-service}}(t)$ describe the
change in $M(t)$ in the time interval $[0,t]$ due to events of new
arrivals, routing, starting service and being in service,
respectively. Hence,
$$M_{\text{arrival}}(t)\geq 0\qquad \mbox{and}\qquad
M_{\text{b-service}}(t)\leq 0.$$ From our definition of the weights, one can
further observe that
\begin{eqnarray}
E[M_{\text{routing}}(t)~|~\cF(0)]~=~0&&\qquad
E[M_{\text{s-service}}(t)~|~\cF(0)]~=~0\notag\\
E[M_{\text{arrival}}(t)~|~\cF(0)]&\leq& t\sum_{d=D}^{\infty}
d\, m_{\max}\,\|P^d\balpha\|_1,\label{eq00:lemkey2}
\end{eqnarray}
where we define $m_{\max}:=\max_{i\in\I} m_i$. 

By appropriately defining constants $\gamma, T, D,
\gamma_1, \gamma_2\in (0,\infty)$, we first prove the following.
\begin{equation}\label{eq0:lemkey2}
E\left[\mathcal G(Y(T))-\mathcal G(Y(0))~|~\cF(0)\right]~\leq~
\begin{cases}
	-\gamma&\mbox{if}~\bQ_{>D}(0)\neq 0\\
	~~~~~~\cC &\mbox{otherwise}
	\end{cases},
\end{equation}
for some constant $\cC<\infty$.
It is not hard to prove \eqref{eq0:lemkey2}
for $\bQ_{>D}(0)=0$,
and hence we only provide the proof for $\bQ_{>D}(0)\neq 0$.
Consider the two complementary events: $\|\bV(0)\|_1\geq 2J B_{\text{renewal}}$ and
$\|\bV(0)\|_1< 2J B_{\text{renewal}}$,
where $B_{\text{renewal}}$ is some constant which will be determined later.

\paragraph{First case.}
On the event $\{\|\bV(0)\|_1\geq 2J B_{\text{renewal}}\}$, we observe that
$\|\bV(t)\|_1=\sum_{j\in\J}V^j(t)$ and
\begin{equation}\label{eq0-1:lemkey2}
E\left[V^j(T)~|~\cF(0)\right]~\leq ~\begin{cases}V^j(0)-\beta_j
T&\mbox{if}~V^j(0)>
\beta_jT\\
\qquad B_{\text{renewal}}
&\mbox{otherwise}\end{cases},\end{equation}
where one can find an appropriate constant $B_{\text{renewal}}<\infty$ depending
on the variances of the generic service times $\{\Gamma_i\}$
using the renewal theory (e.g.,
$B_{\text{renewal}}=\max_{i\in\I} E[\Gamma_i^2]/E[\Gamma_i]$
from the proof of Proposition 6.2 in \cite{asmussen2003applied}).
Hence,
in case there is no $j\in \J$ satisfying $V^j(0)> \beta_j T$,
\begin{eqnarray}
E\left[\|\bV(T)\|_1~|~\cF(0)\right]&=&\sum_{j\in\J} E\left[V^j(T)~|~\cF(0)\right]\notag\\
&\leq&
J B_{\text{renewal}}\notag\\
&\leq&\|\bV(0)\|_1- J B_{\text{renewal}},\label{eq1:lemkey2}
\end{eqnarray}
since $\|\bV(0)\|_1\geq 2J B_{\text{renewal}}$. On the other hand, if there
exists a $j_0\in \J$ satisfying $V^{j_0}(0)> \beta_{j_0} T$, we have
\begin{eqnarray}
E\left[\|\bV(T)\|_1~|~\cF(0)\right]&=&\sum_{j\in\J}
E\left[V^j(T)~|~\cF(0)\right]\notag\\
&=&E\left[V^{j_0}(T)~|~\cF(0)\right]+
\sum_{j\in\J\setminus\{j_0\}}
E\left[V^j(T)~|~\cF(0)\right]\notag\\
&\leq& V^{j_0}(0)-\beta_{j_0}T+ \sum_{j\in\J\setminus\{j_0\}} \max\left\{V^j(0)-\beta_j T,B_{\text{renewal}}\right\}\notag\\
&\leq& V^{j_0}(0)-\beta_{j_0}T+ \sum_{j\in\J\setminus\{j_0\}} \left(V^j(0)+B_{\text{renewal}}\right)\notag\\
&\leq&\|\bV(0)\|_1+ J B_{\text{renewal}}-\beta_{\min}T\notag\\
&\leq&\|\bV(0)\|_1- J B_{\text{renewal}},\label{eq2:lemkey2}
\end{eqnarray}
where we use \eqref{eq0-1:lemkey2} and choose
\begin{equation}\label{eq:defT}
	T~>~\left\lceil \frac{2J B_{\text{renewal}}}{\beta_{\min}} +1\right\rceil.
	\end{equation}
Therefore, in both \eqref{eq1:lemkey2} and \eqref{eq2:lemkey2}, we
have
$$E\left[\|\bV(T)\|_1~|~\cF(0)\right]~\leq~\|\bV(0)\|_1-
J B_{\text{renewal}}.
$$ Using this, it follows that
\begin{eqnarray*}
E[\mathcal G(Y(T))-\mathcal G(Y(0))~|~\cF(0)]&\leq&
E\left[M(T)-M(0)~|~\cF(0)\right]+\gamma_1
E\left[\|\bV(T)\|_1-\|\bV(0)\|_1~|~\cF(0)\right]\\
&\leq&E\left[M_{\text{arrival}}(T)~|~\cF(0)\right]-\gamma_1 J B_{\text{renewal}}\\
&\leq&T\sum_{d=D}^{\infty} d\, m_{\max}\,\|P^d\balpha\|_1
-\gamma_1J B_{\text{renewal}}\\
&\leq&-\gamma,
\end{eqnarray*}
where we use \eqref{eq00:lemkey2} and now define 
$D,\gamma_1$ in terms of $\gamma$ as follows:
\begin{eqnarray*}
D=D(\gamma):=\min\left\{x\in\mathbb{N}~:~T\sum_{d=x}^{\infty}
d\,m_{\max}\,\|P^d\balpha\|_1~\leq~\gamma\right\}\qquad
\gamma_1=\gamma_1(\gamma):=\frac{2\gamma}{J B_{\text{renewal}}}.
\end{eqnarray*}
We specify the value of $\gamma$ at a later stage in the proof.
This completes the proof of \eqref{eq0:lemkey2} on the first event
$\{\|\bV(0)\|_1\geq
2J B_{\text{renewal}}\}$. 

\paragraph{Second case.} Now consider the second event $\{\|\bV(0)\|_1<
{2J B_{\text{renewal}}}\}$. From our choice of $T$ in \eqref{eq:defT}, we have
$$\max_{j\in\J} \frac{V^j(0)}{\beta_j}~\leq~ \frac{\|\bV(0)\|_1}{\beta_{\min}}~<~ \frac{2J
B_{\text{renewal}}}{\beta_{\min}}~<~T-1,$$ which implies that
before time $T-1$, all activities have
completed serving the jobs they were serving at time $0$ (and they could
have worked on other jobs as well). Thus, it follows that
\begin{eqnarray}\label{eq3:lemkey2}
E\left[\|\bV(T)\|_1~|~\cF(0)\right]~\leq~ J B_{\text{renewal}}.\end{eqnarray}
Consider a job with the largest counter
(and hence, contributing to $\bQ_{>D}(0)$) at time $0$ and let $\I^{(h_0)}\subset \I$
be the component the job belongs to. Define the event $\mathcal E_{\text{step 3-1}}$ that
the $\varepsilon$-LRFS policy executes step 3-1 (see the description of $\varepsilon$-LRFS
in Definition \ref{def:epLRFS}) for this component at least once before time $T-1$.
Let $\mathcal E_{\text{step 3-1}}^*$ be the subevent that the job identified in step 1 of the policy
is selected for processing when step 3-1 is carried out for the first time, so that
$\Pr[\mathcal E_{\text{step 3-1}}^*~|~\mathcal E_{\text{step 3-1}}]=\varepsilon$.
On the event $\mathcal E_{\text{step 3-1}}^*$, let
the random variable $j^*$ denote the activity
which is chosen to process this job, and
let $X^*$ be the associated service time. We then have that
$E[X^*|j^*]=m_{i_{j^*}}.$
Observe that
\begin{eqnarray*}
&&E[M_{\text{s-service}}(T)+M_{\text{b-service}}(T)~|~j^*,X^*,\mathcal E_{\text{step 3-1}}^*,\cF(0)]\\
&&\qquad\qquad\qquad\qquad\qquad\qquad\leq~\begin{cases}
-m_{i_{j^*}}+X^*-\beta_{j^*} &\mbox{on the event}~ \{X^*\geq \beta_{j^*}\}\\
\quad\qquad- m_{i_{j^*}} &\mbox{otherwise}
\end{cases}.
\end{eqnarray*}
It thus follows that
\begin{eqnarray*}
\lefteqn{
	E[M_{\text{s-service}}(T)+M_{\text{b-service}}(T)~\big|~\mathcal E_{\text{step 3-1}}^*,\cF(0)]}\\
&\leq&E\left[-m_{i_{j^*}}+[X^*-\beta_{j^*}]_+~\big|~\mathcal E_{\text{step 3-1}}^*,\cF(0)\right]\\
&\leq&E\left[-m_{i_{j^*}}+[X^*-\beta_{\min}]_+~\big|~\mathcal E_{\text{step 3-1}}^*,\cF(0)\right]\\
&=&E\left[E\left[\left.-m_{i_{j^*}}+[X^*-\beta_{\min}]_+~\big|~j^*,X^*,\mathcal E_{\text{step 3-1}}^*,\cF(0)\right]~\right|~\mathcal E_{\text{step 3-1}}^*,\cF(0)\right]\\
&\leq&\max_{i\in \I} E\left[-m_i+[\Gamma_i-\beta_{\min}]_+\right],
\end{eqnarray*}
where $\Gamma_i$ is a generic service time for buffer $i$
and $[x]_+=x$ if $x\geq 0$ and $[x]_+=0$ otherwise.
It is easy to see that
$E\left[-m_i+[\Gamma_i-\beta_{\min}]_+\right]<0$ for all $i$ since $E[\Gamma_i]=m_i>0$ and $\beta_{\min}>0$.
Hence, the above inequality implies that
\begin{equation}\label{eq4-1:lemkey2}
	E[M_{\text{s-service}}(T)+M_{\text{b-service}}(T)~\big|~\mathcal E_{\text{step 3-1}},\cF(0)]
	~\leq~-\varepsilon \nu, \end{equation}
	with $$\nu:=\min_{i\in \I} E\left[m_i-[\Gamma_i-\beta_{\min}]_+\right]>0.$$

If the network is synchronized or if {the second part of Condition {\em C2}} of Theorem \ref{thm:main} holds,
the only way for the event $\mathcal E_{\text{step 3-1}}$ not to occur is that
there exists a processor (for the component $\I^{(h_0)}$) processing a job constantly
during the entire time interval $[\mathcal T^{(h_0)}(0),T-1]$ (i.e., the job starts service
before time $\mathcal T^{(h_0)}(0)$ and is still in service at time $T-1$).
Recall that before time $T-1$, every processor
completes the service requirement of the job it was processing at time $0$.
Hence, in addition to (\ref{eq:defT}), if we choose $T$ to also satisfy
\begin{equation*}
	T~>~2\max_{i\in\I} m_i +2,
\end{equation*}
then it follows that
\begin{equation}\label{eq:event3-1}
	\Pr[\mathcal E_{\text{step 3-1}}~|~\cF(0)]\geq \left(\frac12\right)^{K},
\end{equation}
where we use the fact that if a processor starts to process a new job in the
time interval $(0,\mathcal T^{(h_0)}(0)]\subset(0,1]$,
its service requirement is at most $2\max_{i\in \I}m_i$ with probability $1/2$
by the Markov inequality.
From \eqref{eq4-1:lemkey2} and \eqref{eq:event3-1}, we conclude that
\begin{eqnarray}\label{eq4:lemkey2}
E[M_{\text{s-service}}(T)+M_{\text{b-service}}(T)~|~\cF(0)]~\leq~-\frac{\varepsilon \nu}{2^K}.
\end{eqnarray}

From \eqref{eq00:lemkey2}, \eqref{eq3:lemkey2} and
\eqref{eq4:lemkey2}, the desired inequality \eqref{eq0:lemkey2}
follows as
\begin{eqnarray*}
&&E[\mathcal G(Y(T))-\mathcal G(Y(0))~|~\cF(0)]\\
&&\qquad\leq~
E\left[M(T)-M(0)~|~\cF(0)\right]+\gamma_1
E\left[\|\bV(T)\|_1-\|\bV(0)\|_1~|~\cF(0)\right]\\
&&\qquad\leq~
E\left[M_{\text{arrival}}(T)~|~\cF(0)\right]+
E[M_{\text{s-service}}(T)+M_{\text{b-service}}(T)~|~\cF(0)]+\gamma_1J B_{\text{renewal}}\\
&&\qquad\leq~
T\sum_{d=D}^{\infty} d\,m_{\max}\,\|P^d\balpha\|_1
-\frac{\varepsilon \nu}{2^K}+\gamma_1J B_{\text{renewal}}\\
&&\qquad\leq~
-\gamma,
\end{eqnarray*}
where for the last inequality we define
$$\gamma~:=~\frac{\varepsilon\, \nu}{2^{K+2}}.$$
Here we note that $D=D(\gamma)<\infty$ as long as
$\gamma>0$ (i.e., $\varepsilon>0$).
This completes the proof of
\eqref{eq0:lemkey2} on the second event $\{\|\bV(0)\|_1<
2J B_{\text{renewal}}\}$.

\paragraph{Completing the proof of Lemma \ref{lem:key2}.}
Now \eqref{eq0:lemkey2} implies that
\begin{eqnarray*}
E[\mathcal G(Y(T))^2-\mathcal G(Y(0))^2~|~\cF(0)]&\leq& 2\mathcal
G(Y(0))\,E[\mathcal
G(Y(T))-\mathcal G(Y(0))~|~\cF(0)]+\cC\\
&\leq& 2\left(-\gamma {\bf 1}^+_{\|\bQ_{>D}(0)\|_1}+
\cC \left(1-{\bf 1}^+_{\|\bQ_{>D}(0)\|_1}\right)\right)
\mathcal G(Y(0))+\cC\\
&=& 2\left(-(\gamma+\cC) {\bf 1}^+_{\|\bQ_{>D}(0)\|_1}+
\cC \right)
\mathcal G(Y(0))+\cC
\end{eqnarray*}
where one can check $E\left[\left(\mathcal G(Y(T))-\mathcal
G(Y(0))\right)^2~\big|~\cF(0)\right]<\infty$ for the first inequality and
the precise value of $\cC$ can be different from line to line.
Finally, we define
$$\gamma_2:=m_{\min}\gamma$$
and the conclusion of Lemma \ref{lem:key2} follows from
\begin{eqnarray*}
	E[\mathcal G(Y(T))^2-\mathcal G(Y(0))^2~|~\cF(0)]&\leq&
	2\left(-(\gamma+\cC) {\bf 1}^+_{\|\bQ_{>D}(0)\|_1}+
	\cC\right)
	\mathcal G(Y(0))+\cC\\
&\leq& -2\gamma \,m_{\min}\|\bQ_{>D}(0)\|_1+\cC\left(\|\bQ_{\leq
D}(0)\|_1+\|\bV(0)\|_1+1\right)\\ &=&
-\gamma_2\|\bQ_{>D}(0)\|_1+\cC\left(\|\bQ_{\leq
D}(0)\|_1+\|\bV(0)\|_1+1\right),
\end{eqnarray*}
where we
use $\mathcal G(Y(0))\leq \cC\left(\|\bQ_{\leq
D}(0)\|_1+\|\bV(0)\|_1+1\right)$ for the case $\bQ_{>D}(0)=0$ and
$\mathcal G(Y(0))\geq m_{\min}\|\bQ_{>D}(0)\|_1$ for the other
case $\bQ_{>D}(0)\neq 0$ (i.e., $\|\bQ_{>D}(0)\|_1>0$).

\section*{Acknowledgments}
ABD gratefully acknowledges NSF grant EEC-0926308 for financial support.

\end{document}